\documentclass[11pt, a4paper]{elsarticle}
\usepackage[top=30truemm, bottom=30truemm]{geometry}
\usepackage{newtxtext}

\usepackage{palatino} 

\usepackage{amsmath, amssymb, amsfonts, mathpazo} 
\usepackage{bm}
\usepackage{cases}

\usepackage[driverfallback=dvipdfmx]{hyperref}
\hypersetup{colorlinks=true}

\usepackage{graphicx}
\usepackage{natbib}

\biboptions{authoryear}
\usepackage{framed}
{\endMakeFramed}
\usepackage{mdframed}
\usepackage{amsthm} 
\newtheorem{definition}{Definition}
\newtheorem{lemma}{Lemma}
\newtheorem{assumption}{Assumption}
\newtheorem{corollary}{Corollary}
\newtheorem{proposition}{Proposition}

\usepackage{algpseudocode}
\usepackage{algorithm}

\usepackage{multirow,bigdelim}

\usepackage{listliketab}
\usepackage{booktabs}
\usepackage{tabularx}
\newcolumntype{C}{>{\centering\arraybackslash}X}
\newcolumntype{L}{>{\raggedright\arraybackslash}X}
\newcolumntype{R}{>{\raggedleft\arraybackslash}X}
\usepackage{longtable}

\journal{Transportation Research Part B: Methodological}
\begin{document}

	\begin{frontmatter}
		
		\title{\textbf{Markovian Traffic Equilibrium Assignment\\ based on Network Generalized Extreme Value Model}}
		
		\author[a]{Yuki Oyama\corref{cor1}}
		\ead{oyama@shibaura-it.ac.jp}
		\cortext[cor1]{Corresponding author}
		\author[b]{Yusuke Hara\corref{cor1}}
		\ead{hara@bin.t.u-tokyo.ac.jp}
		\author[c]{Takashi Akamatsu\corref{cor1}}
		\ead{akamatsu@plan.civil.tohoku.ac.jp}
		
		\address[a]{Shibaura Institute of Technology, Department of Civil Engineering, 3-7-5 Toyosu, Koto-ku, Tokyo 135-8548, Japan}
		\address[b]{The University of Tokyo, Next Generation Artificial Intelligence Research Center, 7-3-1 Hongo, Bunkyo-ku, Tokyo 113-8654, Japan}
		\address[c]{Tohoku University, Graduate School of Information Sciences, 06 Aramaki-Aoba, Aoba-ku, Sendai, Miyagi 980-8579, Japan}
		
		\begin{abstract}
			\small
			This study establishes Markovian traffic equilibrium assignment based on the network generalized extreme value (NGEV) model, which we call \textit{NGEV equilibrium assignment}. The use of the NGEV model for route choice modeling has recently been proposed, and it enables capturing the path correlation without explicit path enumeration. 
			However, the theoretical properties of the model in traffic assignment have yet to be investigated in the literature, which has limited the practical applicability of the NGEV model in the traffic assignment field. 
			This study addresses the research gap by providing the theoretical developments necessary for the NGEV equilibrium assignment.
			We first show that the NGEV assignment can be formulated and solved under the same path algebra as the traditional Markovian traffic assignment models. Moreover, we present the equivalent optimization formulations to the NGEV equilibrium assignment. The formulations allow us to derive both primal and dual types of efficient solution algorithms. 
			In particular, the dual algorithm is based on the accelerated gradient method that is for the first time applied in the traffic assignment. The numerical experiments showed the excellent convergence and complementary relationship of the proposed primal-dual algorithms. 
		\end{abstract} 
		
		
	\end{frontmatter}
	
	\newtheorem{thm}{Theorem}
	\newtheorem{prop}{Proposition}
	\newtheorem{defi}{Definition}
	
	\newpage
	\tableofcontents
	\newpage
	
	\section{Introduction}\noindent
	The stochastic traffic equilibrium assignment models provide a static description of how the flows circulate in a congested network where the travel costs are subject to random fluctuations. This randomness may arise from the uncertainties of travelers' perception and of unspecified features in the model. The random utility models (RUMs) that take into account this stochastic nature of travel costs have been utilized in traffic assignment for modeling realistic route choice behavior.
	
	Modeling route choices involves two major problems: path correlation and path enumeration. In a large-scale network, there are many paths overlapping each other, resulting in the correlation among path utilities perceived by travelers. The path enumeration is often required to define a path set and explicitly capture the path correlation. However, the actual path set is unknown to the modeler. In large-scale networks, the path set definition significantly impacts prediction results, and enumerating all feasible paths is computationally impossible. 
	To address these problems, one would need an extension of logit-based models often implemented in practice. A promising approach is the network generalized extreme value (NGEV) model of \cite{Daly2006NGEV}. The NGEV model is the most flexible closed-form model in reproducing path correlation to date.
	\cite{Hara2012NGEV, Hara2014NGEV}, \cite{Papola2013NGEV} and \cite{Mai2016NGEV} recently proposed a novel route choice model based on the NGEV model. This \textit{NGEV route choice model} directly utilizes the transportation network structure as a GEV network and captures the underlying correlation among path utilities. More importantly, it does not require path enumeration, which allows for efficient computation.
	
	Although the NGEV route choice model has notable flexibility 
	and closed-form expression, the theoretical properties of the model in traffic assignment remain uninvestigated. \cite{Papola2013NGEV} performed Dial-based network loading based on the NGEV route choice model. However, the traffic \textit{equilibrium} assignment based on the NGEV route choice model has never been shown in the literature. The implementation of the model in traffic equilibrium assignment requires a convergent and efficient assignment algorithm. The method of successive averages (MSA) \citep{Sheffi1982sue} can be applied to solve the equilibrium assignment with any network loading algorithm, but its poor convergence rate limits its applicability in large-scale networks. A potential way to derive an efficient and convergent algorithm is to formulate an equivalent optimization problem to the definition of equilibrium assignment and analyze its theoretical properties.
	
	The objective of this study is to establish a framework of the Markovian traffic equilibrium assignment based on the NGEV route choice model, which we call the \textit{NGEV equilibrium assignment}. More specifically, our contributions are threefold.
	First, we give a unified perspective between the NGEV route choice model, Markovian traffic assignment (MTA) models, and the deterministic shortest path (SP) model from an algebraic point of view. The path algebra provides generalized formulations for these models, which also shows that the NGEV model is compatible with the efficient loading procedure of MTA. The MTA is as computationally efficient as Dial's algorithm and ensures the consistent solution for the traffic equilibrium assignment \citep{Akamatsu1997Entropy}, and this study proposes its NGEV version.
	Second, we present mathematical definitions and formulations of the NGEV equilibrium assignment. We show that the equivalent optimization problem to the NGEV equilibrium assignment can be formulated based on link-based variables and has a generalized entropy term compared to the logit-based equilibrium assignment. The formulation clarifies mathematical properties of the NGEV equilibrium assignment, such as the uniqueness of the solution and the dual formulation. This is the key to accurately and efficiently solving the NGEV equilibrium assignment.
	Third, we propose both primal and dual types of efficient solution algorithms to the NGEV equilibrium assignment, based on the equivalent optimization formulations. 
	The proposed dual algorithm is based on the accelerated gradient methods that have been recently developed for solving large-scale optimization problems in the machine learning field \citep{Nesterov1983, Beck2009fista}. This study is the first to propose this type of solution algorithm in the traffic assignment field. The proposed algorithm is generally applicable to Markovian traffic equilibrium assignment, including the logit-based assignment.
	Numerical experiments demonstrate the efficiency and excellent convergence properties of the proposed algorithms. In particular, in large-scale or congested networks, the proposed algorithms show incomparable efficiency to the MSA. Moreover, the complementary relationship between the primal and dual algorithms is clarified. The primal algorithm is more efficient when the network is moderately congested, but its performance declines with an increase in the demand level. In contrast, the dual algorithm is not affected by the demand level. When the demand level is high, the dual algorithm is more efficient and has better convergence than the primal algorithm. Moreover, the high computational efficiency of the dual algorithm in large-scale congested networks with many origin-destination (OD) pairs is revealed.
	
	It should be noted that this is a theoretical paper aiming at efficiently solving a highly accurate solution of the NGEV equilibrium assignment. It is not an objective to quickly find a solution with low accuracy. Instead, we are more interested in analyzing the theoretical solution properties and comparing the computational performances of the solution algorithms. Also, even though we use the NGEV model as the underlying route choice model because it has excellent properties of capturing the path correlation and avoiding path enumeration, discussing route choice models (e.g., which model better describes realistic behavior) is not the main focus of this study. Therefore, most of the experiments showed in this paper focus on computational performances of the solution algorithms of the NGEV equilibrium assignment. 
	
	The remainder of this paper is structured as follows. Section \ref{sec:review} provides literature review to emphasize our contributions. Section \ref{sec:routechoice} introduces the NGEV route choice model with a simple illustrative example that shows the model behavior compared to the logit and probit models. Section \ref{sec:loading} proposes an algebra that clarifies the relation between the NGEV assignment and the MTA models. It is also shown that these models can be formulated and solved in a unified manner. Section \ref{sec:formulation} provides the definition and formulations for the NGEV equilibrium assignment, and Section \ref{sec:alg} presents its solution algorithms for both primal and dual formulations. The computational performance of the proposed solution algorithms are examined through numerical experiments in Section \ref{sec:experiments}. In the end, Section \ref{sec:conclusion} concludes the paper. \ref{app:notation} presents the list of notations frequently used in this paper. \ref{sec:proof} provides all the proofs associated with the propositions, corollaries and lemma presented in this paper. 
	\ref{app:loadingtime} - \ref{app:probit} provide a number of supplemental but detailed experiments.
	
	\section{Literature Review}\label{sec:review}
	\subsection{Traffic Assignment Models}\noindent
	Capturing the correlation among path utilities has been one of the main issues of traffic assignment. \cite{Daganzo1977sue} proposed the multinomial probit assignment that represents the correlation by utilizing the covariance matrix. \cite{Yai1997} later proposed a probit model with a structured covariance matrix in order to deal with the complicated path set in railway networks. The probit-based traffic equilibrium assignment was also proposed \citep{Daganzo1979probit, Daganzo1982unconstrained}. 
	However, it is not expressed in a closed-form. 
	To approximately compute the expected link flows, the probit assignment requires the Monte-Carlo simulation that iteratively performs the deterministic assignment as many as it draws link costs from a multivariate normal distribution. This procedure is repeated for the traffic equilibrium assignment. The probit-based assignment models are thus computationally almost intractable to obtain a good approximation in large-scale networks, and this paper hereinafter focuses on the assignment models with a closed-form expression.\footnote{Although limiting the number of draws can reduce the computational effort, it leads to a significant approximation error in large-scale networks where the number of feasible paths between an origin-destination pair is uncountable. Moreover, for the probit-based equilibrium assignment, the MSA is currently the only available solution algorithm, whose poor convergence has been proved in the literature. Therefore, from the computational perspective, obtaining a good approximation of the probit-based equilibrium assignment is almost impossible in large-scale networks. See \ref{app:probit} for more details.}
	
	To retain an efficient computation, a closed-form expression in evaluating route choice probabilities and the use of implicit path enumeration are preferable. \cite{Dial1971} proposed an algorithm for performing logit-based assignment without explicit path enumeration, which has contributed considerably to subsequent studies because of its efficiency. Dial's algorithm restricts the path set to the set of \textit{efficient paths} that never include any moves away from the destination in terms of travel time.
	
	\cite{Bell1995MCA} and \cite{Akamatsu1996MCA} proposed Markovian traffic assignment (MTA) that is also consistent with the logit-based assignment and implicitly considers all feasible alternatives including even cyclic paths. \cite{Akamatsu1997Entropy} proposed the link-based (Markovian) formulation of stochastic equilibrium assignment and showed that the use of MTA models ensures the convergence to an exact equilibrium solution whereas Dial's algorithm does not. \cite{Baillon2008MCA} and \cite{Cominetti2015mte} presented a further investigation into the Markovian traffic equilibrium assignment.
	
	In the context of discrete choice analysis, \cite{Fosgerau2013RL} and \cite{Mai2015NRL} formulated and estimated MTA models, based on the dynamic discrete choice model framework of \cite{Rust1987}. 
	Their models are called ``recursive logit'' (RL) models, and a comprehensive tutorial on RL models is provided in \cite{Zimmermann2020tutorial}. It is worth noting that, although the RL models are derived in a different manner, they are mathematically equivalent to the original MTA models \citep{Akamatsu1996MCA, Baillon2008MCA} that have been studied for over two decades in the traffic assignment field.
	
	The logit-based assignment is unable to capture the underlying correlation among path utilities due to the independence of irrelevant alternatives (IIA) property. 
	A number of route choice models have been proposed based on the generalized extreme value (GEV) model \citep{McFadden1978} to capture the underlying correlation among path utilities in a closed-form expression, such as the cross-nested logit (CNL) model \citep{Vovsha1998CNL} and the paired combinatorial logit (PCL) model \citep{Chu1989paired}. The combination of random utility models (CoRUM) by \cite{Papola2018} is a recent contribution. We have also seen the use of the weibit model in traffic equilibrium assignment \citep{Kitthamkesorn2013W, Nakayama2015W}. 
	However, as investigated by \cite{Prashker1998}, most of such models require explicit path enumeration, which is computationally intractable and inapplicable in practice. 
	
	Another approach to capture the overlap effect is to modify the deterministic utility term of a logit-based assignment. The C-logit model \citep{Cascetta1996Clogit} and the path-size logit (PSL) model \citep{Ben-Akiva1999} introduced correction attributes called a communality factor and a path size attribute respectively, as a function of the physical overlaps with all the other paths in the path set. \cite{Fosgerau2013RL} introduced a link-based correction attribute that uses the expected link flow as a proxy for overlaps. These attributes correct the path utilities to capture the overlap effect, but the models still retain the IIA property.
	
	Recently, \cite{Hara2012NGEV, Hara2014NGEV} and \cite{Papola2013NGEV} proposed a flexible route choice model based on the NGEV model of \cite{Bierlaire2002NGEV} and \cite{Daly2006NGEV}. They assumed a transportation network structure as a direct representation of GEV network, and route choice behavior is modeled as an ordered joint node/link choice from the origin to the destination. The recursive formulation of the NGEV generating function enables capturing the complex correlation structure with a closed-form expression, and moreover, it does not require explicit path enumeration.
	In other words, the NGEV route choice model achieves a balance between efficient and advanced modeling, which has long been the main problem in route choice modeling.
	\cite{Mai2016NGEV} formulated a similar model derived from the dynamic discrete choice framework of \cite{Rust1987}. To represent complex correlation structures, \cite{Mai2016NGEV} considers an NGEV model at each link choice phase by adding artificial states. All the GEV networks are then integrated into a transportation network, and a path choice behavior is modeled on the integrated network. This model can be considered as a more general version, but such network manipulation is costly and unsuitable for traffic assignment in large networks, as compared to the NGEV route choice model by \cite{Hara2012NGEV, Hara2014NGEV} and \cite{Papola2013NGEV} who directly translate a transportation network into a GEV network. In fact, the main application presented in \cite{Mai2016NGEV} is a multilevel CNL formulation that is almost equivalent to the NGEV route choice model.
	
	Despite its notable flexibility, theoretical properties of the NGEV route choice model in the traffic assignment has never been investigated in the literature. Though \cite{Papola2013NGEV} computed Dial-like network loading based on the NGEV route choice model, they do not show the traffic \textit{equilibrium} assignment.. The implementation of the model in traffic equilibrium assignment requires convergent and efficient assignment algorithms, which are discussed in more detail in the next subsection.
	Also, Dial-like algorithms do not ensure convergence to a consistent solution in traffic equilibrium assignment, because they restrict the path set to efficient paths, which vary every time link travel times are updated during iterations of solving the equilibrium assignment \citep{Akamatsu1997Entropy}\footnote{Fixing the set of efficient paths by using constant values such as free-flow travel times is a way for the traffic equilibrium assignment with a Dial-like algorithm to ensure convergence, but in this case a crude and inconsistent solution is obtained in which there would exist paths that have lower travel times but are never loaded.}. 
	The MTA \citep{Akamatsu1996MCA, Baillon2008MCA}, which is also an efficient network loading algorithm but implicitly considers a fixed path set, is mathematically convergent and preferable in terms of the solution properties. Although the MTA has the potential to load excessive flows on cyclic structures of network, it may be alleviated by using the NGEV route choice model. The MTA models have also an obvious similarity to the NGEV route choice model in the formulations \citep{Hara2012NGEV, Hara2014NGEV, Mai2015NRL, Mai2016NGEV, Oyama2019prism}, but the algebraic relation between them has never been investigated and the MTA based on the NGEV route choice has yet to be presented. 
	Moreover, the NGEV equilibrium assignment and its solution algorithms remain uninvestigated. A theoretical development including formulations and solution algorithms for the NGEV equilibrium assignment is required to open up the applicability of NGEV-based traffic assignment models. 
	
	\subsection{Solution Algorithms for Stochastic Traffic Equilibrium Assignment}\noindent
	
	The stochastic traffic equilibrium assignment has long been studied, and a variety of solution algorithms have been proposed. Our review mainly focuses on link-based algorithms that find a solution in the space of link-based variable\footnote{Path-based solution algorithms to the stochastic equilibrium assignment have also been studied in the literature \citep[e.g.,][]{Chen1991algorithms, Damberg1996algorithm}. However, they are not practically applicable in large-scale networks due to the need for path enumeration, which is computationally intractable.}. 
	The first algorithm applied to solve the stochastic equilibrium assignment was the method of successive averages (MSA) \citep{Sheffi1982sue}.
	The MSA can be applied with any stochastic loading algorithm and has been widely used in many different studies. Recently, \cite{Liu2009msa} proposed a modification of MSA using a self-regulated averaging scheme. However, the convergence speed of MSA is slow due to the use of a predetermined sequence of step sizes, which limits its practical applicability. This poor convergence rate of MSA can be improved by optimizing the step size at each iteration. \cite{Chen1991algorithms} proposed an algorithm in which the step size is optimized based on \cite{Fisk1980sue}'s convex minimization problem, but it is suboptimal because only a part of \cite{Fisk1980sue}'s objective function was minimized to avoid evaluating the path-based entropy function. \cite{Maher1997probit} and \cite{Maher1998sue} directly utilized the unconstrained formulation of \cite{Sheffi1982sue} and developed link-based algorithms respectively for probit- and logit-based equilibrium assignment models. 
	Their algorithms determine the optimal step size by using the quadratic interpolation, with the same descent direction as MSA. 
	Moreover, \cite{Akamatsu1997Entropy} proposed the link-based formulation of the logit equilibrium assignment by decomposing the traditional path-based entropy function used in \cite{Fisk1980sue}. Based on this formulation, Akamatsu presented the partial linearization (PL) method \citep{Evans1976PL, Patriksson1993} to solve the stochastic equilibrium assignment. It was proved that the subproblem in PL reduces to performing logit-based network loading, and that the step size is efficiently optimized because the whole objective function of \cite{Akamatsu1997Entropy} can be evaluated only by link-based variables. \cite{Lee2010pl} presented a modification of the PL method based on the descent direction improvement strategy of \cite{Fukushima1984fw} and the parallel tangent (PARTAN) technique \citep{Arezki1990fw}. The link-based PL method, however, has never been applied to the previously mentioned GEV-based assignment models that can address the path correlation problem because such models are often presented with the path-based formulation \citep[e.g.][]{Bekhor1999formulations}. 
	
	Although the dual formulation that finds the equilibrated link costs was early presented by \cite{Daganzo1982unconstrained}, most of the existing solution algorithms to stochastic equilibrium assignment can be categorized as primal algorithms in the sense that their unknown is in the space of flow variable.
	\cite{Maher1997probit}, \cite{Maher1998sue}, and \cite{Xie2012dual} directly utilized \cite{Sheffi1982sue}'s unconstrained model, which is a flow-based alternative to the dual formulation by \cite{Daganzo1982unconstrained}, but they still find the solution in the flow space. There are, though, several possible advantages of directly solving the dual formulation. The dual problem reduces to unconstrained maximization programing, and its gradient can be efficiently evaluated by stochastic loading. Also, the dimension of solution space of the cost variable is smaller than that of the flow variable.
	These advantages inspired us to solve the stochastic equilibrium assignment by using the dual-type gradient-based methods. Gradient-based methods for large-scale optimization problems have been recently studied and have contributed significantly to the rapid growth of research in the machine learning field. In the stochastic traffic assignment context, only an application of the primitive gradient projection (GP) method can be found in \cite{Bekhor2005review}. They applied the GP method to a relaxed problem of \cite{Fisk1980sue}'s model, but their solution algorithm was based on the path flow variable and thus required path enumeration, which was computationally expensive. Despite their potential, dual-type gradient-based methods have not been well investigated in the traffic assignment field. In particular,  accelerated gradient methods \citep[e.g.,][]{Nesterov1983, Beck2009fista, Su2014differential, Donoghue2015adaptive} that outperform GP have never been used but can potentially provide an efficient solution algorithm to the stochastic traffic equilibrium assignment.
	
	\section{NGEV Route Choice Model}\label{sec:routechoice}\noindent
	This section introduces the NGEV route choice model, which captures the underlying path correlation without explicitly enumerating path alternatives. 
	We transform a transportation network structure into a GEV network \citep{Daly2006NGEV}, upon which the NGEV route choice model is formulated. A simple illustrative example to show the model property is also provided.
	
	\subsection{Defining GEV Network}\noindent
	This study directly uses the topological structure of a transportation network to define a GEV network in order to capture the underlying path correlation.
	Let $\mathcal G \equiv ({\mathcal N}, {\mathcal L})$ be a digraph representing a transportation network where ${\mathcal N}$ and ${\mathcal L}$ are the sets of nodes and links respectively. We define the set of successor nodes $\mathcal{F}(i)$ and the set of predecessor nodes $\mathcal{B}(i)$ for each node $i \in {\mathcal N}$. If link $ij$ exists in ${\mathcal L}$, then node $i$ is connected to node $j$, resulting in $j \in \mathcal{F}(i)$ and $i \in \mathcal{B}(j)$.
	Each link $ij \in \mathcal{L}$ is associated with the generalized link travel cost $c_{ij}$, which can be flow-dependent. 
	Moreover, let ${\mathcal O} \subseteq {\mathcal N}$ denote the set of origin nodes, and ${\mathcal D} \subseteq {\mathcal N}$ denote the set of destination nodes. The set of OD pairs is defined as ${\mathcal W}$.
	
\begin{figure}[htbp]
		\begin{center}
			\includegraphics[width=12cm]{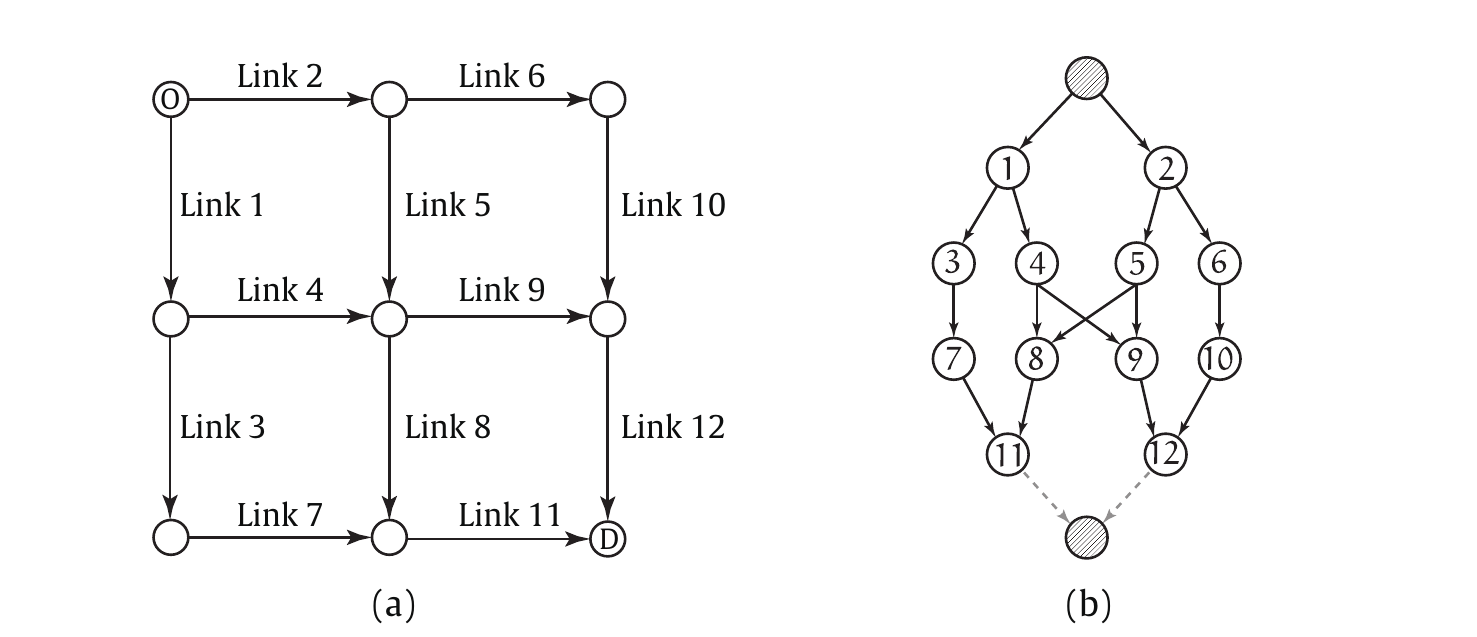}
			\caption{(a) A transportation network , and (b) its transformation into a GEV network where the links of the transportation network are represented as the vertexes.}
			\label{fig:illustrative_example}
		\end{center}
	\end{figure}

	Given an OD pair $(o,d)$$\in {\mathcal W}$, a route choice on network $\mathcal{G}$ can be described as a joint choice of all elemental links of the route in their ordered sequence from the origin $o$ towards the destination $d$, which \cite{Papola2013NGEV} call an \textit{ordered joint choice context}. In such a context, the structure of transportation network $\mathcal{G}$ corresponds to a GEV network as follows: the origin $o$ corresponds to the \textit{root}, which is faced with the initial choice among the outgoing links $\{oj, \forall j \in \mathcal{F}(o)\}$ from $o$ ; each link $ij \in \mathcal{L}$ corresponds to a \textit{vertex}, which is faced with a choice among the outgoing links $\{jk, \forall k \in \mathcal{F}(j)\}$ from node $j$ ; and the links heading to the destination node $d$ correspond to the final vertexes on a GEV network. 
	Therefore, any ordered sequence of vertexes connecting the root to a final vertex identifies a particular path of the transportation network originating from $o$ and terminating at $d$. Fig. \ref{fig:illustrative_example} shows an example of transformation of a transportation network structure into a GEV network. 
	Note that any transportation network structure can be transformed into a GEV network in the same way.
	
	\subsection{Route Choice Model}\noindent
	Following \cite{Papola2013NGEV}, we formulate a route choice model on a GEV network based on an ordered joint choice context through multilevel cross-nested structures. 
	Each choice of link $ij$ conditional on node $i$ represents an elemental choice of the route, associated with the elemental utility $U^d_{ij} = - (c_{ij} + \varepsilon^d_{ij})$ where the error term $\varepsilon^d_{ij}$ follows a multivariate extreme value distribution with the generating function $G^d_{ij}$ defined by the NGEV model. 
	Given an OD pair $(o, d) \in {\mathcal W}$, it is assumed that a traveler would choose, among the set $\mathcal{R}^{od}$ of routes available between $(o, d)$, the route minimizing the generalized travel cost. Under the link-additive cost assumption, the generalized travel cost is given by the sum of the elemental utilities of all the elemental links $ij \in \mathcal{L}_r$ of $r$:
	\begin{equation}
	c_r + \varepsilon_{r} = \sum_{ij \in \mathcal{L}_r} (c_{ij} + \varepsilon^d_{ij}) = \sum_{ij \in \mathcal{L}_r} c_{ij} + \sum_{ij \in \mathcal{L}_r} \varepsilon^d_{ij}.
	\label{eq:cr}
	\end{equation}
	Moreover, in the joint choice context, the choice probability $p(r)$ of route $r$ is given by the product of the conditional probabilities $p^d_{ij|i}$ of link $ij$ on node $i$:
	\begin{equation}
	p(r) = \prod_{ij \in \mathcal{L}_r} p^d_{ij|i},
	\label{eq:pathprob}
	\end{equation}
	and the NGEV model defines $p^d_{ij|i}$ as \citep{Daly2006NGEV}:
	\begin{equation}
	p^d_{ij|i} = \frac{\alpha^d_{ji} (G^d_{ij})^{\theta^d_i/\theta^d_j}}
	{{\displaystyle \sum_{j' \in \mathcal{F}(i)}} \alpha^d_{j'i} (G^d_{ij'})^{\theta^d_i/\theta^d_{j'}}},
	\label{eq:pijG}
	\end{equation}
	where $\theta^d_i$ is the variance scale parameter associated with node $i$, which is strictly positive, and $\alpha^d_{ji}$ is the allocation parameter that represents the degree of membership of node $j$ to predecessor node $i \in \mathcal{B}(j)$, satisfying $\sum_{i \in \mathcal{B}(j)} \alpha^d_{ji} = 1$, $\alpha^d_{ji} \ge 0$. In addition, the NGEV generating function $G^d_{ij}$ associated with vertex (link) $ij$ is recursively formulated as follows:
	\begin{equation}
	G^d_{ij} = 
	\begin{cases}
	e^{- \theta^d_j c_{ij}} {\displaystyle \sum_{k \in \mathcal{F}(j)}} \alpha^d_{kj} (G^d_{jk})^{\theta^d_j/\theta^d_k} & {\rm if} \:\: j \neq d,\\
	e^{- \theta^d_j c_{ij}} & {\rm if} \:\: j = d.\\
	\end{cases}
	\label{eq:G}
	\end{equation}
	
	The approach of this study is to bridge between the NGEV route choice model (\ref{eq:pathprob})-(\ref{eq:G}) and the link-based formulation of MTA \citep{Akamatsu1997Entropy}. 
	To this end, we here translate the NGEV generating function $G^d_{ij}$ of link $ij$ into the \textit{expected minimum cost} $\mu^d_j$ of node $j$ to destination $d$ in the following equation:
	\begin{equation}
	G^d_{ij} = e^{\theta^d_j \bar{U}^d_{ij}} = e^{- \theta^d_j (c_{ij} + \mu^d_j)},
	\label{eq:Gmu}
	\end{equation}
	where $\bar{U}^d_{ij} = - (c_{ij} + \mu^d_j)$ is the expected maximum utility associated with link $ij$ \citep{Daly2006NGEV}.
	Finally, from (\ref{eq:pijG})-(\ref{eq:Gmu}), the NGEV route choice model is given as :
	\begin{align}
	&p^d_{ij|i} = \frac{\alpha^d_{ji} e^{- \theta^d_i (c_{ij} + \mu^d_j)}}
	{\sum_{j' \in \mathcal{F}(i)} \alpha^d_{j'i} e^{- \theta^d_i (c_{ij'} + \mu^d_{j'})}}
	&\forall ij \in \mathcal{L}, \ \forall d \in \mathcal{D},
	\label{eq:pij}\\
	&1 = \sum_{j \in \mathcal{F}(i)} \alpha^d_{ji} e^{- \theta^d_i (c_{ij} + \mu^d_j - \mu^d_{i})}
	&\forall i \in \mathcal{N}, \ \forall d \in \mathcal{D},
	\label{eq:mu}
	\end{align}
	whose solution is a set $\bm{\mu} \equiv \{\bm{\mu}^d, \forall d \in \mathcal{D}\}$, where $\bm{\mu}^{d} \equiv (\mu^d_i)_{i \in \mathcal{N}}$ is the vector of the expected minimum costs from nodes to destination $d$. We will discuss its relation and give a unified perspective with the existing MTA models in Section \ref{sec:algebra}.
	
	The NGEV route choice model can describe an arbitrary correlation structure by assuming different values of parameters $\theta^d_i$ and $\alpha^d_{ji}$. Moreover, it does not require explicit path enumeration because a route choice probability is evaluated as the product of link choice probabilities of (\ref{eq:pij}), as shown in (\ref{eq:pathprob}). These explain the notable flexibility and efficiency of the NGEV route choice model. Note that it is easily proved that the NGEV route choice model (\ref{eq:pij}) and (\ref{eq:mu}) corresponds to the logit-based route choice model as a special case, when assuming that $\theta^d_i = \theta, \forall i \in \mathcal{N}$ and $\alpha^d_{ji} = 1, \forall ij \in \mathcal{L}$. It also reduces to the SP model when $\theta^d_i \rightarrow \infty, \forall i \in \mathcal{N}$ and $\alpha^d_{ji} = 1, \forall ij \in \mathcal{L}$.
	
	We refer the reader to \cite{Papola2013NGEV} for detailed experiments to analyze the characteristics of an NGEV route choice model. They tested the capability of an NGEV route choice model reproducing the effects of path overlaps in four different networks and as compared to a number of other route choice models including the logit, C-logit, PSL, and CNL models as well as the probit model.

	\section{Network Loading and Path Algebra}\label{sec:loading}\noindent
	This section discusses the network loading, i.e. \textit{flow-independent} assignment, based on the NGEV route choice model, which we call the \textit{NGEV assignment} (we later define a \textit{flow-dependent} version as the \textit{NGEV equilibrium assignment}).
	To begin, we define the mathematical conditions of the NGEV assignment based on the NGEV route choice formulations presented in Section \ref{sec:routechoice}. We then introduce a path algebra that provides a unified perspective to the NGEV assignment with existing MTA models. An illustrative example to compare some loading models is also presented.
	
	\subsection{NGEV Assignment}\noindent
	Let $\mathbold{q}^d \equiv (q^d_i)_{i \in \mathcal{N}}$ be the given demand flows toward destination $d \in \mathcal{D}$: $q^d_i > 0$ when $(i,d) \in \mathcal{W}$, and $q^d_i = 0$ otherwise. 
	The NGEV assignment assigns the given demand $\{\mathbold{q}^d, \forall d \in \mathcal{D}\}$ on the network and determines a set $\mathbold{x} \equiv \{\mathbold{x}^d, \forall d \in \mathcal{D}\}$ where $\mathbold{x}^d \equiv (x^d_{ij})_{ij \in \mathcal{L}}$ is the vector of the destination $d$-specific link flows. 
	At every node $i$, the expected inflow $z^d_i$ that node $i$ receives and that is directed toward $d$ is distributed among the successor links $\{ij, j \in \mathcal{F}(i)\}$ according to the link choice probabilities $p^d_{ij|i}$ of the NGEV route choice model (\ref{eq:pij}) and (\ref{eq:mu}):
	\begin{align}
	x^d_{ij} =  p^d_{ij|i} z^d_i.
	\label{eq:pflow}
	\end{align}
	This is a link-based and many-to-one traffic assignment procedure. Then the node flow $z^d_i$ is given by the sum of expected link flows $x^d_{hi}$ incoming from nodes $h \in \mathcal{B}(i)$ and directed toward $d$, plus the demand flow $q^d_{i}$, which is expressed as the following flow conservation equation:
	\begin{align}
	z^d_{i} = q^d_i + \sum_{h \in \mathcal{B}(i)} x^d_{hi}.
	\label{eq:z}
	\end{align}
	From (\ref{eq:pflow}) and (\ref{eq:z}), for any node $i \in \mathcal{N}$ and destination $d \in \mathcal{D}$, the following flow conservation law should hold:
	\begin{equation}
	\sum_{j \in \mathcal{F}(i)} x^d_{ij} - \sum_{h \in \mathcal{B}(i)} x^d_{hi}
	= 
	\begin{cases}
	q^d_i & {\rm if} \:\: i \neq d,\\
	- \sum_{o \in \mathcal{O}}q^d_o & {\rm if} \:\: i = d.\\
	\end{cases}
	\label{eq:stateflow}
	\end{equation}
	In addition, for any link $ij \in \mathcal{L}$ and destination $d \in \mathcal{D}$, the link flow $x_{ij}^d$ should satisfy the usual non-negativity constraint:
	\begin{equation}
	x^d_{ij} \ge 0.
	\label{eq:edgeflow}
	\end{equation}
	More concisely, Eqs.(\ref{eq:stateflow}) and (\ref{eq:edgeflow}) can be written in the following vector forms, for all $d \in \mathcal{D}$:
	\begin{align}
	&\mathbf{A}\mathbold{x}^d = \tilde{\mathbold{q}}^d,
	\label{eq:stateflow_vector} \\
	& \mathbold{x}^d \ge 0,
	\label{eq:edgeflow_vector}
	\end{align}
	where $\mathbf{A} \equiv [a_{k,ij}]_{k \in \mathcal{N}, ij \in \mathcal{L}}$ is the node-link incidence matrix: $a_{k,ij} = 1$ if link $ij$ is outgoing from node $k$, $a_{k,ij} = -1$ if link $ij$ is incoming to node $k$, and $a_{k,ij} = 0$ otherwise. 
	$\tilde{\mathbold{q}}^d \equiv (\tilde{q}^d_i)_{i \in \mathcal{N}}$ is the vector of modified demand flows: $\tilde{q}^d_i = q^d_i$ for $i \neq d$, and $\tilde{q}^d_d = -\sum_{o \in \mathcal{O}}q^d_o$.
	Finally, the definition of the (flow-independent) NGEV assignment condition is given as follows:
	\begin{definition}{\textbf{(NGEV assignment)}}\label{def:sta}
		The NGEV assignment is mathematically defined as a problem finding a solution tuple $\langle \mathbold{x}, \bm \mu \rangle$ that satisfies (\ref{eq:pij})-(\ref{eq:edgeflow}).
	\end{definition}
	
	As the loading procedure of the NGEV assignment, the link flows $\mathbold{x}$ can be calculated based on the link choice probability matrix $\mathbf{P} \equiv [p_{ij|i}]_{i,j \in \mathcal{N}}$ by (\ref{eq:pij}). 
	Substituting (\ref{eq:pflow}) into (\ref{eq:z}) yields:
	\begin{equation}
	\mathbold{z}^d = \mathbf{P}^{\top} \mathbold{z}^d + \mathbold{q}^d,
	\label{eq:zmatrix}
	\end{equation}
	which is a system of linear equations whose solution is the vector of node flows $\mathbold{z}^d$ \citep{Akamatsu1996MCA, Baillon2008MCA}. Once $\mathbold{z}^d$ is solved, then the link flow $\mathbold{x}^d$ can be computed using (\ref{eq:pflow}). 
	
	\subsection{Path-algebra for Markovian Traffic Assignment}\label{sec:algebra}\noindent
	To evaluate the link choice probability $\mathbf{P}$ of (\ref{eq:pij}), the expected minimum cost $\bm \mu$ has to be solved.
	As expressed in (\ref{eq:mu}), the expected minimum cost is solved through a recursive formulation, meaning that the NGEV assignment is based on the dynamic programming (DP) framework. This implies the relation between the NGEV assignment and the MTA models.
	To further discuss this relation, we enumerate below the formulations of the (expected) minimum costs $\bm \mu_{\rm{SP}}$, $\bm \mu_{\rm{Logit}}$, and $\bm \mu_{\rm{NGEV}}$, respectively of the SP assignment \citep{Floyd1962, Warshall1962}, the logit assignment \citep{Akamatsu1996MCA, Fosgerau2013RL}, and the NGEV assignment:
	\begin{align}
	&\mu^d_{i, \rm{SP}} = \min_{j \in \mathcal{F}(i)} \{ c_{ij} + \mu^d_{j, \rm{SP}} \},
	\label{eq:SPDP} \\
	&\mu^d_{i, \rm{Logit}} = \mathbb{E}\left[\min_{j \in \mathcal{F}(i)} \{ c_{ij}  + \mu^d_{j, \rm{Logit}} + \varepsilon_{ij} \}\right]
	= -\frac{1}{\theta} \ln \sum_{j \in \mathcal{F}(i)} e^{-\theta (c_{ij} + \mu^d_{j, \rm{Logit}})},
	\label{eq:MCADP} \\
	&\mu^d_{i, \rm{NGEV}} 
	= \mathbb{E}\left[\min_{j \in \mathcal{F}(i)} \{ c_{ij} - \frac{1}{\theta^d_i} \ln \alpha^d_{ji}  + \mu^d_{j, \rm{NGEV}} + \varepsilon_{ij} \}\right]
	=  
	-\frac{1}{\theta^d_i} \ln \sum_{j \in \mathcal{F}(i)} \alpha^d_{ji} e^{-\theta^d_i (c_{ij} + \mu^d_{j, \rm{NGEV}})}.
	\label{eq:NGEVDP}
	\end{align}
	Note that (\ref{eq:NGEVDP}) is a transformation of (\ref{eq:mu}). As shown above, every MTA model is formulated in the form of recurrence relation and has the same structure of a fixed point problem: $\bm \mu = F(\bm \mu)$. Furthermore, it is worth mentioning that (\ref{eq:SPDP})-(\ref{eq:NGEVDP}) all share a common algebraic structure: \textit{the mapping $F: \mathbb{R}^{|\mathcal{N}|} \rightarrow \mathbb{R}^{|\mathcal{N}|}$ for each model can be regarded as ``linear'' from an abstract algebraic point of view.}
	
	To show this point in more detail, we introduce an extended version of the path algebra by \cite{Carre1979graphs}. 
	We define the algebra $\mathbb{R}_{\rm path} \equiv \langle \mathbb{R}_{\epsilon}, \oplus, \otimes \rangle$ as a set $\mathbb{R}_{\epsilon}= \mathbb{R} \cup \{\epsilon\}$ equipped with two binary operations $\oplus$ and $\otimes$, which have the following elementary properties: 
	\begin{align*}
	&\textbf{(a) commutative law:} 
	\hspace{1.0cm}
	x \oplus y = y \oplus x   
	\hspace{1.0cm}
	x \otimes y = y \otimes x  
	\hspace{2.5cm}
	\forall x,y \in \mathbb{R}_{\epsilon} \\
	&\textbf{(b) associative law:} 
	\hspace{0.3cm}
	(x \oplus y) \oplus z = x \oplus (y \oplus z) 
	\hspace{0.3cm}
	(x \otimes y) \otimes z = x \otimes (y \otimes z) 
	\hspace{0.3cm}
	\forall x,y,z \in \mathbb{R}_{\epsilon}\\
	&\textbf{(c) distributive law:} 
	\hspace{2.0cm}
	x \otimes (y \oplus z) = (x \otimes y) \oplus (y \otimes z) 
	\hspace{2.1cm}
	\forall x,y,z \in \mathbb{R}_{\epsilon}
	\end{align*}
	Moreover, set $\mathbb{R}_{\epsilon}$ contains a \textit{zero element} $\epsilon$ such that $x \oplus \epsilon = \epsilon \oplus x = x$ and $x \otimes \epsilon = \epsilon \otimes x = \epsilon$, $\forall x \in \mathbb{R}_{\epsilon}$, and a \textit{unit element} $e$ such that $x \otimes e = e \otimes x = x$, $\forall x \in \mathbb{R}_{\epsilon}$. An \textit{inverse} $x^{\otimes -1} = -x$ that satisfies $x \otimes x^{\otimes -1} = x^{\otimes -1} \otimes x = e$ exists for any $x \in \mathbb{R}$. 
	In other words, algebra $\mathbb{R}_{\rm path}$ is a \textit{semiring} in terms of algebraic structure.
	
	For this algebra, we may also define matrix operations. Let $\mathcal{M}_{n}(\mathbb{R}_{\epsilon})$ be the set of all $n \times n$ matrices whose entries belong to $\mathbb{R}_{\epsilon}$. The two binary operations on $\mathcal{M}_{n}(\mathbb{R}_{\epsilon})$ are defined as follows:
	\begin{equation}
	\mathbf{A \oplus B} \equiv [a_{ij} \oplus b_{ij}]
	\hspace{0.5cm}
	\mathrm{and}
	\hspace{0.5cm}
	\mathbf{A \otimes B} \equiv \left[\bigoplus^n_{k=1} a_{ik} \otimes b_{kj}\right];
	\label{eq:matrix_operation}
	\nonumber
	\end{equation}
	and the powers of a matrix are: 
	\begin{equation}
	\mathbf{A}^{\otimes m} \equiv \mathbf{A}^{\otimes m-1} \otimes \mathbf{A}
	\hspace{0.5cm}
	\mathrm{and}
	\hspace{0.5cm}
	\mathbf{A}^{\otimes 0} \equiv \mathbf{E},
	\nonumber
	\label{eq:matrix_power}
	\end{equation}
	where $\mathbf{E}$ is the \textit{unit matrix} of $\mathcal{M}_{n}(\mathbb{R}_{\epsilon})$, whose element $e_{ij}$ equals the unit element $e$ if $i = j$ and the zero element $\epsilon$ otherwise. $\mathbf{E}$ also satisfies $\mathbf{E} \otimes \mathbf{A} = \mathbf{A} \otimes \mathbf{E} = \mathbf{A}$ for any matrix $\mathbf{A} \in \mathcal{M}_{n}(\mathbb{R}_{\epsilon})$.
	
	Having defined algebra $\mathbb{R}_{\rm path}$, we now state that every model of (\ref{eq:SPDP})-(\ref{eq:NGEVDP}) can be written in a unified manner:
	\begin{equation}
	\label{eq:mu_path}
	\bm \mu^d = \mathbf{W} \otimes \bm \mu^d \oplus \mathbold{e}^d
	\end{equation}
	where $\mathbf{W} \equiv [w_{ij}] \in \mathcal{M}_{|\mathcal{N}|}(\mathbb{R}_{\epsilon})$ is the link weight matrix, and $\mathbold{e}^d$ is the column vector, corresponding to node $d$, of $\mathbf{E}$. Eq.(\ref{eq:mu_path}) is consistent with each model of (\ref{eq:SPDP})-(\ref{eq:NGEVDP}) when defining the set $\mathbb{R}_{\epsilon}$, operations $\oplus, \otimes$, and elements of the link weight matrix according to Table \ref{tab:path_choice}\footnote{
		We herein define node-specific operations $\oplus_i, \otimes_i, \forall i \in \mathcal{N}$, so as to describe the node-specific scale parameter of the NGEV assignment. The corresponding matrix operations are:
		\begin{equation} \nonumber
		\mathbf{A \oplus B} \equiv [a_{ij} \oplus_i b_{ij}]
		\hspace{0.5cm}
		\mathrm{and}
		\hspace{0.5cm}
		\mathbf{A \otimes B} \equiv \left[\bigoplus^n_{k=1} \!_i  a_{ik} \otimes_i b_{kj}\right].
		\end{equation}
		For the SP and logit cases, we simply assume $\oplus_i \equiv \oplus, \otimes_i \equiv \otimes, \forall i \in \mathcal{N}$.
	}.  
	Operation $\oplus$ for the NGEV assignment can be interpreted as a generalization of the min function $\min \{x, y\}$ for the SP assignment so that it returns the expected minimum cost.
	
	\begingroup
	\renewcommand{\arraystretch}{1.5}
	\begin{table}[h]
		\caption{Definitions of algebra $\mathbb{R}_{\rm path}$ for MTA models.}
		\label{tab:path_choice}
		\footnotesize
		\begin{tabular*}{\hsize}{@{\extracolsep{\fill}}|c|c|c|c|@{}}
			\hline
			\textbf{Definition}&\hspace{-6pt}\textbf{Shortest path} &\hspace{-6pt}\textbf{Logit} &\textbf{Network GEV} \\
			\hline
			$\mathbb{R}_{\epsilon}$ & $\mathbb{R}\cup \{\infty\}$ & $\mathbb{R}\cup \{\infty\}$ & $\mathbb{R}\cup \{\infty\}$ \\
			\hline
			$x\oplus_i y$ & $\min \{x,y \}$ & $ -\frac{1}{\theta} \ln [\exp(-\theta x) + \exp(-\theta y)]  $ &  $  -\frac{1}{\theta_i} \ln [\exp(-\theta_i x) + \exp(-\theta_i y)] $
			\\
			\hline
			$x\otimes_i y$ &  $x + y$ & $-\frac{1}{\theta} \ln [\exp(-\theta x) \cdot \exp(-\theta y)]$ & $-\frac{1}{\theta_i} \ln [\exp(-\theta_i x) \cdot \exp(-\theta_i y)]$ 
			\\
			\hline
			$w_{ij}$ & 
			$\left\{ \begin{array}{ll} c_{ij} & ij \in \mathcal{L} \\ \infty &  ij \not\in \mathcal{L}\\ \end{array} \right.$ & 
			$\left\{ \begin{array}{ll} c_{ij} & ij \in \mathcal{L} \\ \infty &  ij \not\in \mathcal{L}\\ \end{array} \right.$ & 
			$\left\{ \begin{array}{ll} c_{ij} - \frac{1}{\theta_i} \ln \alpha_{ji}& ij \in \mathcal{L} \\ \infty &  ij \not\in \mathcal{L}\\ \end{array} \right.$ 
			\\
			\hline
		\end{tabular*}
	\end{table}
	\endgroup
	
	Furthermore, we observe a certain similarity between (\ref{eq:mu_path}) and \textit{a system of linear equations} $\bm \mu = \mathbf{W} \bm \mu + \mathbold{e}$ in ordinary matrix algebra. It naturally follows that, regardless of which model we use, Eq.(\ref{eq:mu_path}) can be solved in an analogous manner to that for the ordinary system of linear equations: if the sequence of powers of $\mathbf{W}$ is convergent, by recursively substituting (\ref{eq:mu_path}) into its right-hand side , we have (herein we omit the superscript $d$ for $\bm \mu$ and $\mathbold{e}$ as Eq.(\ref{eq:mu_path}) has the same structure for all $d \in \mathcal{D}$):
	\begin{align}
	\label{eq:V}
	\bm \mu 
	&= \mathbold{e} \oplus \mathbf{W} \otimes \bm \mu
	= \mathbold{e} \oplus \mathbf{W} \otimes [\mathbold{e} \oplus \mathbf{W} \otimes \bm \mu] \nonumber\\
	&= [\mathbf{E} \oplus \mathbf{W}] \otimes \mathbold{e} \oplus \mathbf W^{\otimes 2} \otimes \bm \mu
	= [\mathbf{E} \oplus \mathbf{W}] \otimes \mathbold{e} \oplus \mathbf W^{\otimes 2} \otimes [\mathbold{e} \oplus \mathbf{W} \otimes \bm \mu] \nonumber\\
	& = [\mathbf{E} \oplus \mathbf{W} \oplus \mathbf W^{\otimes 2}] \otimes \mathbold{e} \oplus \mathbf W^{\otimes 3} \otimes \bm \mu \nonumber\\
	& = \cdots \nonumber \\
	& = [\mathbf{E} \oplus \mathbf{W} \oplus \mathbf W^{\otimes 2} \oplus \mathbf W^{\otimes 3} \oplus \cdots] \otimes \mathbold{e}.
	\end{align} 
	That is, solving (\ref{eq:mu_path}) reduces to calculating a matrix power series of $\mathbf{W}$ under the algebra $\mathbb{R}_{\rm path}$, which yields the (expected) minimum cost $\bm \mu$. 
	Once $\bm \mu$ is obtained, we can compute the link choice probability matrix $\mathbf{P}$, and finally the link flows $\mathbold x$ through (\ref{eq:zmatrix}) and (\ref{eq:pflow}). 
	
	Note that, as long as the operations and link weight matrix are appropriately defined, any MTA models \citep[e.g.][]{Mai2015NRL, Oyama2017GRL} other than SP, logit, and NGEV can be solved in the same manner, and Eq.(\ref{eq:V}) always returns the corresponding expected minimum cost for each model. With respect to other various traditional path problems that have the same algebraic structure, see \cite{Carre1979graphs} and \cite{Baras2010path}; in addition, see \cite{Baccelli1992} and \cite{Heidergott2014maxplus} for the Max-Plus algebra, which is a particular instance of the path algebra. 
	
	
	\subsection{Illustrative Example}\label{sec:loading_example}\noindent
	Finally, we provide an example to analyze characteristics of the NGEV assignment with MTA that is presented in the previous subsection, by comparing two assignment models, logit and NGEV, respectively with two different loading algorithms, MTA and Dial's algorithm \citep{Dial1971} (which we refer to as `Dial' for simplicity). In general, Dial and MTA implicitly assume different sets of paths, respectively efficient paths and the universal set, and we know each of them has a drawback. Dial assigns flows only on efficient paths, which is a common problem of logit \citep{Dial1971} and NGEV assignment \citep{Papola2013NGEV}. MTA may assign excessive flows on cyclic paths, but the NGEV formulation can alleviate the problem. Although this experiment is a limited case, it shows a potential of the NGEV assignment with MTA addressing the limitations of logit assignment with MTA as well as of Dial's algorithm.
	
	Several different parameter settings are tested for NGEV assignments. To summarize, the following four cases are tested (with both Dial and MTA, i.e., eight scenarios in total):
	\begin{itemize}
		\item \textbf{Model 1}: Logit assignment, i.e., $\theta_i \equiv 1, \forall i \in \mathcal{N}$, $\alpha_{ji} \equiv 1, \forall ij \in \mathcal{L}$;
		\item \textbf{Model 2}: NGEV assignment with $\theta_i \equiv \frac{D^d(o)}{D^d(i)}, \forall i \in \mathcal{N}$, $\alpha_{ji} \equiv \frac{1}{|\mathcal{B}(j)|}, \forall ij \in \mathcal{L}$;
		\item \textbf{Model 3}: NGEV assignment with $\theta_i \equiv \frac{\pi}{\sqrt{3 D^d(i)}}, \forall i \in \mathcal{N}$, $\alpha_{ji} \equiv \frac{1}{|\mathcal{B}(j)|}, \forall ij \in \mathcal{L}$;
		\item \textbf{Model 4}: NGEV assignment with $\theta_i \equiv \frac{\pi}{\sqrt{6 D^d(i)}}, \forall i \in \mathcal{N}$, $\alpha_{ji} \equiv \frac{1}{|\mathcal{B}(j)|}, \forall ij \in \mathcal{L}$,
	\end{itemize}
	where $D^d(i)$ is the SP cost from node $i$ to destination $d$, upon which Dial's efficient paths are also defined.
	The scale parameter settings for Models 3 and 4 follow \cite{Papola2013NGEV}. 
	
\begin{figure}[htbp]
		\begin{center}
			\includegraphics[width=12cm]{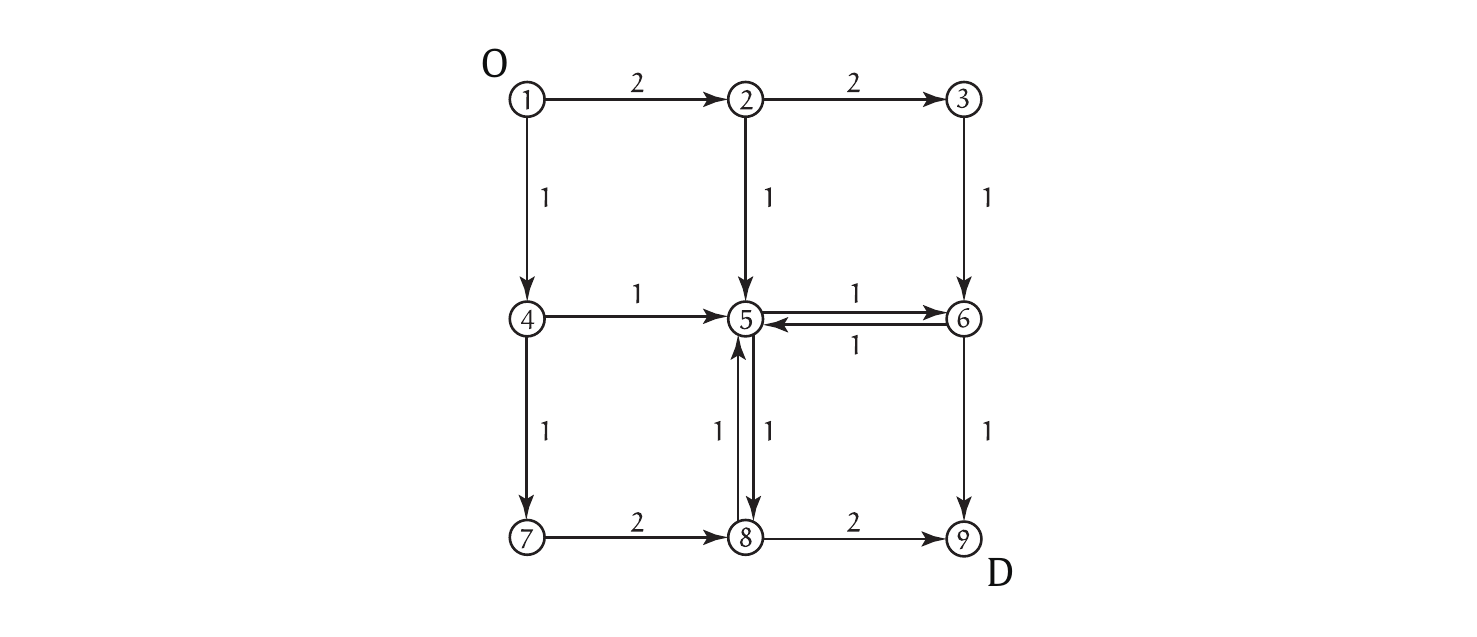}
			\caption{A simple cyclic network where the numbers on links denote link costs.}
			\label{fig:loading_network} 
		\end{center}
	\end{figure}

	Fig. \ref{fig:loading_network} shows a simple cyclic network used for the example, where the number associated with each link on the network indicates the link cost $c_{ij}$. We assume a unit demand flow $q^d_o = 1.0$ for a single OD pair $(1,9)$. Table \ref{tab:loading_examples} reports the loading results of the eight scenarios, indicating the computed link flows (or flow rates as the OD flow equals $1.0$). 
	
	As known, Dial and MTA assume different path sets, namely, the set of efficient paths and the universal set, and therefore, their results are different from each other. When Dial is used as the loading algorithm, links 4-7, 7-8, 5-8, and 8-9 are regarded \textit{inefficient} and have zero flows, which is unreasonable and in particular problematic for the equilibrium assignment.\footnote{Those links should have some flows in stochastic traffic assignment, in particular in the flow-dependent case. When fixing efficient paths with reference to the free-flow travel costs, the equilibrium assignment never assigns flows on those links. Therefore, to assign flows on the links in the equilibrium assignment, the efficient paths have to be defined based on flow-dependent link costs. However, it causes the fluctuation of the path set in the solution process and does not ensure the convergence.}
	This problem of Dial remains in the NGEV models, because it is caused by the path set definition.
	The MTA provides a solution to this limitation by implicitly taking all feasible paths into account. 
	However, in the logit case (Model 1), the MTA excessively assigns flows on routes involving overlaps, and even on cyclic paths. This causes large flows on links 2-5, 4-5, 5-6, 5-8, 8-5, and 6-5 (Model 1 with MTA). 
	The NGEV assignment with MTA may address the respective problems of Dial and logit-MTA. As shown in Table \ref{tab:loading_examples}, Models 2-4 with MTA assign some flows on links 4-7, 7-8, 5-8, and 8-9, and few flows on links 8-5 and 6-5, which is more reasonable. This result indicates that MTA, unlike Dial, is complemented by the NGEV formulation.
	
	Models 2 and 3 show similar behavior with each other. Their behavior are more deterministic than that of Model 4 as their scales $\bm \theta$ take larger values. This is also why Model 4 with MTA assigns small but non-negligible flows on links 8-5 and 6-5. Model 3, though, is more efficient because it uses destination-specific scale parameters, whereas Model 2 uses OD-specific parameters. In fact, with the destination-specific parameters, NGEV assignment requires only as much computational effort as logit assignment (though this fact is clear from the discussion in Section \ref{sec:algebra}, we refer the reader for more details to Appendix \ref{app:loadingtime} that reports an experiment to compare the computational times of network loading). 
	
	It is worth noting that the MTA may become computationally intractable when a network includes cyclic structures, as reported in \cite{Oyama2017GRL, Oyama2019prism} with some illustrative examples. It is also possible that the MTA loads excessive flows on cyclic paths in a network including many cycles, though capturing the correlation by the NGEV model can somewhat alleviate it. Nevertheless, the fact that the path set is fixed regardless of travelers or link costs is mathematically preferable and ensures the convergence of equilibrium assignment to a consistent solution \citep{Akamatsu1997Entropy}.
	
	\begin{table}[htbp]
		\centering 
		\footnotesize
		\caption{Loading results of eight different scenarios.}
		\label{tab:loading_examples}
		\begin{tabular}{| c | c | c | c | c | c |c | c |c|c|c|c|c|c|c|c|}\hline
			&&\multicolumn{13}{c}{\textbf{Link flow rates}}&\\\hline
			\textbf{Model} & \textbf{Loading} &1-2&1-4&2-3&2-5&3-6&4-5&5-6&4-7&5-8&6-9&7-8&8-9&8-5&6-5 \\
			\hline
			1 (Logit) & MTA &.29&.71&.06&.23&.06&.63&.83&.08&.41&.71&.08&.29&\textbf{.20}&\textbf{.18} \\
			1 (Logit) & Dial
			&.33&.67&.09&.24&.09&.67&.91&\textbf{.00}&\textbf{.00}&\textbf{1.0}&\textbf{.00}&\textbf{.00}&.00&.00 \\
			\hline
			2 (NGEV) & MTA
			&.36&.64&.18&.18&.18&.57&.68&.08&.07&.86&.08&.14&\textbf{.00}&\textbf{.00} \\
			2 (NGEV) & Dial
			&.39&.61&.20&.19&.20&.61&.80&\textbf{.00}&\textbf{.00}&\textbf{1.0}&\textbf{.00}&\textbf{.00}&.00&.00 \\
			\hline 
			3 (NGEV) & MTA
			&.39&.61&.21&.18&.21&.47&.54&.14&.13&.75&.14&.25&\textbf{.01}&\textbf{.00} \\
			3 (NGEV) & Dial
			&.46&.54&.27&.19&.27&.54&.73&\textbf{.00}&\textbf{.00}&\textbf{1.0}&\textbf{.00}&\textbf{.00}&.00&.00 \\
			\hline 
			4 (NGEV) & MTA
			&.45&.55&.27&.18&.27&.37&.43&.19&.15&.69&.19&.31&\textbf{.03}&\textbf{.01} \\ 
			4 (NGEV) & Dial
			&.57&.43&.37&.20&.37&.43&.63&\textbf{.00}&\textbf{.00}&\textbf{1.0}&\textbf{.00}&\textbf{.00}&.00&.00 \\
			\hline
		\end{tabular}
	\end{table}
	
	\section{NGEV Equilibrium Assignment}\label{sec:formulation}\noindent
	This section presents the formulations of the \textit{NGEV equilibrium assignment}. We first present the definition of the NGEV equilibrium assignment, and then show its equivalent optimization problem. The Lagrangian dual formulation is finally derived.
	
	\subsection{Assumption and Definition}\noindent
	The NGEV equilibrium assignment is defined based upon the \textit{flow-independent} NGEV assignment conditions (\textbf{Definition \ref{def:sta})}.
	We herewith define the aggregated link flows $\mathbold{X} \equiv (X_{ij})_{ij \in \mathcal{L}}$ which satisfies:
	\begin{equation}
	\mathbold{X} = \sum_{d \in \mathcal{D}} \mathbold{x}^d.
	\label{eq:linkflow_vector}
	\end{equation}
	An additional assumption on the \textit{flow-dependent} link cost function is introduced:
	\begin{assumption}
		\label{assumption:cX}
		The generalized link cost $\mathbold{c}$ is given as a function of link flows $\mathbold{X}$\footnote{$\mathbold{c(X)}$ is a generalized description and includes the case in which cost $c_{ij}$ of link $ij$ is influenced only by flow $X_{ij}$ of link $ij$, i.e., $\mathbold{c} = (c_{ij}(X_{ij}))_{ij \in \mathcal{L}}$.}, i.e., $\mathbold{c} \equiv \mathbold{c}(\mathbold{X}): \mathbb{R}_+^{|\mathcal{L}|} \rightarrow \mathbb{R}_+^{|\mathcal{L}|}$, and the Jacobian of $\mathbold{c}(\mathbold{X})$ is symmetric.
		$\mathbold{c}(\mathbold{X})$ is a continuous single-valued, strictly monotone function, which satisfies:
		\begin{equation}
		(\mathbold{c}(\mathbold{X}) - \mathbold{c}(\mathbold{Y})) \cdot (\mathbold{X - Y}) > 0
		\hspace{0.5cm}
		\forall \mathbold{X \neq Y} \in \mathbb{R}_+^{|\mathcal{L}|}
		\hspace{0.5cm}
		\label{eq:dcdX}
		\end{equation}
	\end{assumption}
	Note that this is a standard assumption in the traffic equilibrium assignment to ensure the uniqueness of the solution, since the seminal works of \cite{Smith1979existence} and \cite{Dafermos1980}.
	The assumption along with the NGEV assignment conditions leads to the definition of NGEV equilibrium assignment:
	\begin{definition}{\textbf{(NGEV equilibrium assignment)}}\label{def:sue}
		The NGEV equilibrium assignment is mathematically defined as a problem finding a solution tuple $\langle \mathbold{c, x}, \bm \mu \rangle$ that satisfies (\ref{eq:pij})-(\ref{eq:edgeflow}), (\ref{eq:linkflow_vector}) and (\ref{eq:dcdX}).
	\end{definition}
	
	\subsection{Equivalent Optimization Problem}\noindent
	Having defined the NGEV equilibrium assignment, we present its equivalent optimization problem. As a preliminary, we first introduce the \textit{flow-independent} case, i.e., the NGEV assignment, in which the link cost $c_{ij}$ is a constant.
	The equivalent optimization problem to the NGEV assignment (\ref{eq:pij})-(\ref{eq:edgeflow}) is formulated as follows.
	
	\begin{mdframed}
		\vspace{-0.25cm}
		\begin{proposition}
			\label{prop:DTA}
			The NGEV assignment conditions (\ref{eq:pij})-(\ref{eq:edgeflow}) are equivalent to finding $\mathbold{x}$ that solves the following optimization problem \textbf{[NGEV]}:
			\vspace{\baselineskip}
			
			\noindent{\bf [NGEV]}
			\begin{equation}
			\label{eq:ZNGEV}
			\min_{\mathbold{x} \ge 0} Z(\mathbold{x})
			\equiv
			\sum_{d \in \mathcal{D}} 
			[
			\mathbold{c} \cdot \mathbold{x}^d - \hat{\bm{\theta}}^d \cdot \mathbold{H}^d(\mathbold{x}^d)
			]
			\end{equation}
			s.t., for all $d \in \mathcal{D}$:
			\begin{equation}
			\mathbf{A}\mathbold{x}^d = \tilde{\mathbold{q}}^d, \nonumber
			\end{equation}
			where 
			\begin{equation*}
			\hat{\bm{\theta}}^d \equiv \left(\hat{\theta_i}\right)_{i \in \mathcal{N}}, \mathbold{H}^d(\mathbold{x}^d) \equiv 
			\left( H^d_i(\mathbold{x}^d) \right)_{i \in \mathcal{N}}, \mathbold{z}^d \equiv
			\left( z^d_i \right)_{i \in \mathcal{N}}
			\end{equation*}
			with 
			\begin{align}
			&\hat{\theta_i} \equiv \frac{1}{\theta^d_i} \\
			&H^d_i(\mathbold{x}^d) \equiv -\sum_{j \in \mathcal{F}(i)} x^d_{ij}
			\ln \frac{x^d_{ij}}{\alpha^d_{ji} z^d_i}
			\label{eq:entropy} \\
			&z^d_i \equiv \sum_{j \in \mathcal{F}(i)} x^d_{ij}
			\label{eq:zi}
			\end{align}
			
			
		\end{proposition}
	\end{mdframed}
	
	Owing to the above equivalent optimization formulation, we can now show a mathematical property of the flow patterns computed by the NGEV assignment. 
	
	\begin{corollary}
		\label{coro:DTAunique}
		The globally optimal solution for the problem \textbf{[NGEV]} can be uniquely determined.
	\end{corollary}
	
	\textbf{Corollary \ref{coro:DTAunique}} states that the NGEV assignment (\ref{eq:pij})-(\ref{eq:edgeflow}) always yields a unique flow pattern for each cost pattern. Next, a straightforward extension of the above discussion leads to the equivalent optimization problem of the \textit{flow-dependent} case, i.e., the NGEV equilibrium assignment:
	\begin{mdframed}
		\vspace{-0.25cm}
		\begin{proposition}
			\label{prop:SUE}
			The NGEV equilibrium assignment conditions (\ref{eq:pij})-(\ref{eq:edgeflow}), (\ref{eq:linkflow_vector}) and (\ref{eq:dcdX}) are equivalent to finding $\mathbold{x}$ that solves the following optimization problem \textbf{[NGEV-FD/P]}:
			\vspace{\baselineskip}
			
			\noindent{\bf [NGEV-FD/P]}
			\begin{equation}
			\label{eq:NGEVSUE}
			\min_{\mathbold{x} \ge 0} Z^{\rm FD}_{\rm P}(\mathbold{x})
			\equiv
			C(\mathbold{X})
			- 
			\sum_{d \in \mathcal{D}} \hat{\bm{\theta}}^d \cdot \mathbold{H}^d(\mathbold{x}^d),
			\end{equation}
			s.t.,
			\begin{align}
			&\mathbf{A}\mathbold{x}^d = \tilde{\mathbold{q}}^d,
			\nonumber\\
			&\mathbold{X} = \sum_{d \in \mathcal{D}} \mathbold{x}^d, \nonumber
			\end{align}
			where
			\begin{equation}
			C(\mathbold{X}) \equiv 
			\oint_{\mathbold{X}} \mathbold{c(X)} \mathrm{d}\mathbold{X}.
			\label{eq:Cint}
			\end{equation}
		\end{proposition}
	\end{mdframed}
	
	As shown in (\ref{eq:NGEVSUE}), the problem \textbf{[NGEV-FD/P]} is a slight modification of the problem \textbf{[NGEV]}. Only the first term of the objective function (\ref{eq:ZNGEV}) is modified into the integral term (\ref{eq:Cint}), and the second term remains unchanged. 
	Under \textbf{Assumption \ref{assumption:cX}}, a mathematical property of the problem \textbf{[NGEV-FD/P]} with respect to the uniqueness of solution can be clarified, similarly to the flow-independent case.
	
	\begin{corollary}
		\label{coro:SUEunique}
		The globally optimal solution for problem \textbf{[NGEV-FD/P]} can be uniquely determined.
	\end{corollary}
	
	\subsection{Dual Formulation}\noindent
	To further discuss the properties of the NGEV equilibrium assignment, we herewith present the Lagrangian dual formulation of the problem \textbf{[NGEV-FD/P]}.
	From the discussion in \textbf{Proposition \ref{prop:SUE}} and its proof (\ref{app:proof_prp2}), we obtain the following Lagrangian dual problem:
	\begin{equation}
	\max_{\bm{\mu}^d \in \mathcal{K}_d, \forall d \in \mathcal{D}}
	Z^{\rm FD}_{\rm D}(\bm{\mu})
	\equiv \min_{\mathbold{x} \ge 0} L(\mathbold{x}, \bm{\mu})
	\label{eq:Zdual1}
	\end{equation}
	where $\mathcal{K}_d$ is the feasible region of $\bm{\mu}^d$, defined as:
	\begin{equation}
	\mathcal{K}_d \equiv \{ \bm{\mu}^d \in \mathbb{R}^{|\mathcal{N}|} | 1 = \sum_{j \in \mathcal{F}(i)} \alpha^d_{ji} \exp[-\theta^d_i(c_{ij} + \mu^d_j - \mu^d_i)] \hspace{0.2cm} \forall i \in \mathcal{N} \}.
	\label{eq:K}
	\end{equation}
	By defining $\bm{\mu}^d$ as a function of $\mathbold{c}$, i.e., $\bm{\mu}^d \equiv \bm{\mu}^d(\mathbold{c}): \mathbb{R}_+^{|\mathcal{L}| \rightarrow \mathbb{R}^{|\mathcal{N}|}}$, these further reduce to the following maximization problem with the unknown $\mathbold{c}$: 
	\begin{mdframed}
		\vspace{-0.25cm}
		\begin{proposition}
			\label{prop:SUE_dual}
			The dual problem of \textbf{[NGEV-FD/P]} is given by the following maximization problem \textbf{[NGEV-FD/D]}: 
			\vspace{\baselineskip}
			
			\noindent{\bf [NGEV-FD/D]}
			\begin{equation}
			\max_{\mathbold{c}\ge \underline{\mathbold{c}}} 
			Z^{\rm FD}_{\rm D}(\mathbold{c}) 
			\equiv
			- C^*(\mathbold{c})
			+ \sum_{d \in \mathcal{D}} \bm{\mu}^d(\mathbold c) \cdot \tilde{\mathbold{q}}^d, \label{eq:sue/dual}
			\end{equation}
			where $C^*(\mathbold{c})$ is the conjugate dual (Legendre transform) of $C(\mathbold{X})$, defined as
			\begin{equation}
			C^*(\mathbold{c}) \equiv 
			\max_{\mathbold{X}} \left[ \mathbold{c} \cdot \mathbold{X} - C(\mathbold{X}) \right] = \oint_{\mathbold{c}} \mathbold{c^{-1}(\mathbold{c})} \mathrm{d}\mathbold{c},
			\label{eq:cstar}
			\end{equation}
			and $\mathbold{c}^{-1}(\mathbold{c}): \mathbb{R}^{|\mathcal{L}|}_+ \rightarrow \mathbb{R}^{|\mathcal{L}|}_+$ is the inverse of the link cost function $\mathbold{c}(\mathbold{X})$.
			
		\end{proposition}
	\end{mdframed}
	
	Note that the problem \textbf{[NGEV-FD/D]} has a globally optimal and unique solution if assuming an one-to-one correspondence between $\mathbold{X}$ and $\mathbold{c}$.
	Furthermore, we obtain the following lemma concerning the Lagrangian dual problem \textbf{[NGEV-FD/D]}:
	\begin{mdframed}
		\vspace{-0.25cm}
		\begin{lemma}
			\label{lemma:dual}
			The objective function $Z^{\rm FD}_{\rm D}$ is smooth and concave with respect to $\mathbold c$. The gradient $\nabla_{\mathbold c} Z^{\rm FD}_{\rm D}$ is given by
			\begin{equation}
			\nabla_{\mathbold c} Z^{\rm FD}_{\rm D}(\mathbold c)
			= \mathbold{X}\mathbold{(c)} - \mathbold{c}^{-1}\mathbold{(c)}
			\label{eq:gradient}
			\end{equation}
			where $\mathbold{X(t)} \equiv \sum_{d \in \mathcal{D}} \mathbold{x}^d\mathbold{(t)}$ is the (total) link flow pattern obtained from the NGEV assignment based on link cost pattern $\mathbold t$.
		\end{lemma}
	\end{mdframed}

	\textbf{Lemma \ref{lemma:dual}} shows that the gradient of $Z^{\rm FD}_{\rm D}$ represents the difference between the link demand function $\mathbold{X(c)}$ and the inverse link cost function $\mathbold{c}^{-1}\mathbold{(c)}$, which may be interpreted as an excess demand.

	\section{Solution Algorithms}\label{sec:alg}\noindent 
	This section describes solution algorithms for the NGEV equilibrium assignment.
	We first propose an algorithm based on the primal formulation \textbf{[NGEV-FD/P]}: the partial linearization (PL) method. The PL method is a well-known algorithm used to efficiently solve the logit equilibrium assignment, and we show that it can also be applied to the NGEV equilibrium assignment. We then propose another algorithm based on the dual formulation \textbf{[NGEV-FD/D]}. This dual algorithm is based on the accelerated gradient projection (AGP) method, which has never been applied in the traffic assignment field, even to the logit equilibrium assignment.
	As previously discussed in Section \ref{sec:formulation}, the main difference between the primal and dual problems is in their unknowns: the link flow pattern $\mathbold{x}$ for the primal problem; and the link cost pattern $\mathbold{c}$ for the dual problem.
	
	\subsection{Primal Algorithm: Partial Linearization Method}\noindent
	The PL method updates the current solution based on the descent direction vector at each iteration. The descent direction at the $m$-th iteration is determined by solving the following partially linearized subproblem of \textbf{[NGEV-FD/P]}:
	\begin{equation}
	\label{eq:ZPL}
	\min_{\mathbold{y} \ge 0}  
	\left\{
	Z(\mathbold{y}) \equiv
	\sum_{d \in \mathcal{D}}
	[ \mathbold{c}(\mathbold{x}^{(m)}) \cdot \mathbold{y}^d
	- \hat{\bm{\theta}}^d \cdot \mathbold{H}\mathbold{(y)}
	]
	\mid
	\mathbf{A}\mathbold{y}^d = \tilde{\mathbold{q}}^d, \forall d \in \mathcal{D}
	\right\}
	\end{equation}
	where the solution $\mathbold{y}^*$ is the auxiliary link flow pattern, and $\mathbold{x}^{(m)}$ is the current link flow pattern at the $m$-th iteration. The descent direction vector is then given by $\mathbold{d} = \mathbold{y}^* - \mathbold{x}^{(m)}$.
	This subproblem corresponds to the problem \textbf{[NGEV]} based on link cost pattern $\mathbold{c}(\mathbold{x}^{(m)})$; that is, solving (\ref{eq:ZPL}) is equivalent to the NGEV assignment, which is performed by the algorithm proposed in Section \ref{sec:loading}.
	
	Using this, the PL method is summarized as follows:
	\begin{mdframed}
		\textbf{Partial Linearization (PL)}
		\begin{description}
			\item[Step 0: Initialization.]
			Set $m=0$ and $\mathbold{c}^{(0)} = \underline{\mathbold{c}}$. Assign the OD-flows $\{\mathbold{q}^d,  \forall d \in \mathcal{D}\}$ by the NGV assignment based on the initial link cost pattern $\mathbold{c}^{(0)}$ and obtain the initial solution $\mathbold{x}^{(0)}$.
			\item[Step 1: Cost update.]
			Update the link cost based on the current solution $\mathbold{x}^{(m)}$ by
			$\mathbold{c}^{(m)} := \mathbold{c}(\mathbold{x}^{(m)})$.
			\item[Step 2: NGEV assignment for direction finding.] 
			Solve the subproblem (\ref{eq:ZPL}) by performing the NGEV assignment for $\{\mathbold{q}^d,  \forall d \in \mathcal{D}\}$ based on the current link cost pattern $\mathbold{c}^{(m)}$, and obtain the auxiliary flow pattern $\mathbold{y}^{(m)}$. Determine the descent direction $\mathbold{d}^{(m)} = \mathbold{y}^{(m)} - \mathbold{x}^{(m)}$.
			\item[Step 3: Step size determination.]
			Determine the step size $\gamma^*$
			by solving the following line search problem:
			\begin{equation}
			\gamma^* = \arg \min_{0 \le \gamma \le 1} Z^{\rm FD}_{\rm P}(\mathbold{x}^{(m)} + \gamma \mathbold{d}^{(m)}).
			\nonumber
			\end{equation}
			\item[Step 4: Solution update.] 
			Update the solution by $\mathbold{x}^{(m+1)} := \mathbold{x}^{(m)} + \gamma^* \mathbold{d}^{(m)}$
			\item[Step 5: Convergence test.] 
			If the convergence criterion holds, stop. Otherwise, set $m := m + 1$ and return to \textbf{Step 1}.
		\end{description}
	\end{mdframed}
	For convex programming, the convergence of the PL method to a unique and globally optimal solution is guaranteed \citep{Patriksson1993}.
	
	Note that, for the primal problem \textbf{[NGEV-FD/P]}, the MSA is also applicable to the NGEV equilibrium assignment. The MSA defines the step size $\gamma^* \equiv \frac{1}{m+1}$ at the $m$-th iteration instead of solving the line search problem. In general, MSA is not practically applicable due to the poor convergence. In contrast, the efficiency of PL has been reported in the literature \citep{Patriksson1993, Akamatsu1997Entropy, Lee2010pl}. We will later report numerical experiments to show the efficiency of PL compared to MSA for the NGEV equilibrium assignment.
	
	\subsection{Dual Algorithm: Accelerated Gradient Projection Method}\noindent
	Next, we propose an algorithm for efficiently solving the 
	Lagrangian dual problem \textbf{[NGEV-FD/D]}: an accelerated gradient projection (AGP) method. This algorithm is based on a gradient (first-order) method because computing the Hessian of $Z^{\rm FD}_{\rm D}(\mathbold{c})$ in (\ref{eq:sue/dual}) is almost impossible in large-scale networks. It should be noted that even a single evaluation of the objective function $Z^{\rm FD}_{\rm D}(\mathbold{c})$ or computing the gradient $\nabla Z^{\rm FD}_{\rm D}$ requires performing the NGEV assignment (see \textbf{Lemma \ref{lemma:dual}}), which is computationally expensive in large-scale networks.
	
	The AGP method is an application of accelerated proximal gradient methods, which have been developed recently in machine learning field \citep[e.g.,][]{Beck2009fista,Su2014differential,Donoghue2015adaptive}. These accelerated methods are based on the seminal study by \cite{Nesterov1983}, which showed that a minor modification (the choice of step size and the addition of an extra momentum step) of the gradient methods achieves the known complexity bound $\Omega(1/k^2)$, i.e., in the worst case, any iterative method based solely on the function and gradient evaluations cannot achieve a better accuracy than $\Omega(1/k^2)$ at iteration $k$. 
	
	For the NGEV equilibrium assignment, the dual algorithm considers the following projection problem at each iteration: 
	\begin{align}
	\mathbold{c}^{(m+1)} 
	&:= \arg \min_{\mathbold{c} \in \mathcal{C}} \left( -Z^{\rm FD}_{\rm D}(\mathbold{c}^{(m)}) - \nabla_{\mathbold c} Z^{\rm FD}_{\rm D}(\mathbold{c}^{(m)}) \cdot (\mathbold{c} - \mathbold{c}^{(m)}) + \frac{1}{2s} ||\mathbold{c} - \mathbold{c}^{(m)}||^2 \right) \nonumber \\
	&\equiv {\rm Proj}_C \left( \mathbold{c}^{(m)} + s \nabla_{\mathbold c} Z^{\rm FD}_{\rm D}(\mathbold{c}^{(m)}) \right) 
	\label{eq:projection}
	\end{align}
	where $s$ should be chosen so as to satisfy $0 < s \le 1/L$, for the convergence of the algorithm, and $L$ is a Lipschitz constant of $\nabla_{\mathbold c} Z^{\rm FD}_{\rm D}$. From \textbf{Lemma \ref{lemma:dual}}, the gradient $\nabla_{\mathbold c} Z^{\rm FD}_{\rm D}$ of the problem \textbf{[NGEV-FD/D]} can be obtained by performing the NGEV assignment.
	Thus, considering the feasible region $\mathcal{C}$ of $\mathbold{c}$, we can see that the projection operation in (\ref{eq:projection}) reduces to:
	\begin{align}
	\mathbold{c}^{(m+1)} &:= {\rm Proj}_C \left( \mathbold{c}^{(m)} + s (\mathbold{X}(\mathbold{c}^{(m)}) - \mathbold{c}^{-1} (\mathbold{c}^{(m)}) \right) \nonumber \\
	&= \left< \mathbold{c}^{(m)} + s (\mathbold{X}(\mathbold{c}^{(m)}) - \mathbold{c}^{-1} (\mathbold{c}^{(m)})) \right>_{\underline{c}+}
	\label{eq:update}
	\end{align}
	where $ \mathbold{b} = \left<\mathbold{a}\right>_{\underline{c}+}$ denotes the following element-wise operation $b_{ij} = \max \{a_{ij}, \underline{c}_{ij}\}$, $\forall ij \in \mathcal{L}$, and $\underline{c}_{ij}$ denotes the free-flow travel cost of link $ij$. 
	
	Whereas the gradient projection (GP) method simply iterates updating the current solution by (\ref{eq:update}) with an arbitrary chosen value of $s$, 
	the AGP method modifies the updating phase of GP using a momentum term. 
	The AGP method for \textbf{[NGEV-FD/D]} is summarized as follows:
	\begin{mdframed}
		\textbf{Accelerated Gradient Projection (AGP)}
		\begin{description}
			\item[Step 0: Initialization.] Set: $m :=0$, $j :=0$, $\mathbold{c}^{(0)} := \underline{\mathbold{c}}$, $\mathbold{b}^{(0)} := \underline{\mathbold{c}}$, $t_0:=1$, $0 < s \le 1/L$. 
			\item[Step 1: NGEV assignment.] Assign the OD-flows $\{\mathbold{q}^d,  \forall d \in \mathcal{D}\}$ by the NGV assignment based on the current link cost pattern $\mathbold{b}^{(m)}$, which yields a link flow pattern $\mathbold{X}^{(m)} = \sum_{d \in \mathcal{D}} \mathbold{x}^{d(m)}$.
			\item[Step 2: Updating.] Update the current solution using a momentum term.
			\begin{align}
			\mathbold{c}^{(m+1)} &:= {\rm Proj}_c \left( \mathbold{b}^{(m)} + s \nabla_{\mathbold c} Z^{\rm FD}_{\rm D}(\mathbold{b}^{(m)}) \right)
			\nonumber \\
			&= \left<\mathbold{b}^{(m)} + s (\mathbold{X}^{(m)} - \mathbold{c}^{-1} (\mathbold{b}^{(m)}))\right>_{\underline{c}+},\\
			t_{j+1} &= \frac{1+\sqrt{1+4 t_j^2}}{2},\\
			\mathbold{b}^{(m+1)} &= \mathbold{c}^{(m+1)} + \frac{t_j -1}{t_{j+1}} (\mathbold{c}^{(m+1)} - \mathbold{c}^{(m)}).
			\end{align}
			\item[Step 3: Adaptive restart.]
			If the restart criterion hold, $j := 0$. Otherwise, $j := j+1$.
			\item[Step 4: Convergence test.]
			If the convergence criterion holds, stop. Otherwise, set $m := m +1$ and return to \textbf{Step 1}.
		\end{description}
	\end{mdframed}
	A few remarks are in order here. First, as the restart criterion for the adaptive restart, this paper uses the \textit{function scheme}: for $j \ge k_{\min}$,
	\begin{equation}
	Z^{\rm FD}_{\rm D}(\mathbold{c}^{(m+1)}) < Z^{\rm FD}_{\rm D}(\mathbold{c}^{(m)})
	\end{equation}
	where $k_{\min}$ is the minimum number of iterations to restart, i.e., the adaptive restart never occurs if $j$ is smaller than $k_{\min}$. Other schemes are also available as listed in \cite{Donoghue2015adaptive}.
	Second, the step size $s$ can be adjusted at each iteration by using a backtracking procedure described in \cite{Beck2009fista}. Specifically, at the $m$-th iteration, we find the smallest non-negative integer $i_m$ such that, with $s = \xi^{i_m} s_{m-1}$,
	\begin{equation}
	Z^{\rm FD}_{\rm D}(p_{s}(\mathbold{b}^{(m)})) 
	\le Q_s (p_{s}(\mathbold{b}^{(m)}), \mathbold{b}^{(m)})
	\label{eq:z_BT}
	\end{equation}
	where
	\begin{equation}
	p_{s}(\mathbold{b}^{(m)}) \equiv {\rm Proj}_c \left( \mathbold{b}^{(m)} + s \nabla_{\mathbold c} Z^{\rm FD}_{\rm D}(\mathbold{b}^{(m)}) \right),
	\label{eq:p_BT}
	\end{equation}
	\begin{equation}
	Q_s (\mathbold{x}, \mathbold{y})
	\equiv
	Z^{\rm FD}_{\rm D}(\mathbold{y}) 
	+ (\mathbold{x} - \mathbold{y}) \cdot \nabla_{\mathbold c} Z^{\rm FD}_{\rm D}(\mathbold{y})
	+ \frac{1}{2s} ||\mathbold{x} - \mathbold{y}||^2.
	\label{eq:Q_BT}
	\end{equation}
	In numerical experiments in the next section, this backtracking procedure is actually applied, and a comparison of the cases with and without the procedure is provided in \ref{app:sensitivity}.
	
	\section{Numerical Experiments} \label{sec:experiments}\noindent
	This section presents numerical experiments of the NGEV equilibrium assignment to show the computational performances of the algorithms proposed in Section \ref{sec:alg}. Note that the experiments focus on solution algorithms for the NGEV equilibrium assignment, but we refer the reader to a comparison of the flow patterns of the logit and NGEV equilibrium assignment in \ref{app:sue}. 
	
	For all the experiments, we set the parameters of the NGEV route choice model as $\theta_i = \frac{\pi}{\sqrt{3 D^d(i)}}$, $\forall i \in \mathcal{N}$ and $\alpha_{ji} = \frac{1}{|\mathcal{B}(j)|}$, $\forall ij \in \mathcal{L}$, which are the same as Model 3 described in Section \ref{sec:loading_example}, and implemented the MTA as the loading module. 
	The link cost function is defined as $c_{ij}(X_{ij}) = \underline{c}_{ij} [1 + ( \frac{X_{ij}}{\kappa_{ij}})^4]$ for all link $ij \in \mathcal{L}$, where $\underline{c}_{ij}$ is the free-flow link cost and $\kappa_{ij}$ is the nominal link capacity. 
	We solve the line search problem in PL (\textbf{Step 3}) by the golden section method with a threshold of $10^{-3}$. 
	All the algorithms have been implemented in Python 3.6 on a machine with 14 core Intel Xeon W processors (2.5 GHz) and 64 GB of RAM. 
	
	\subsection{Algorithms Convergence}\noindent
	We first discuss the convergence behavior of the proposed algorithms, using
	the Sioux Falls network data provided by \cite{SiouxFalls2016}. The network data contain 24 nodes, 76 links, 576 OD pairs, and 360,600 trips. 
	The free-flow link cost $\underline{c}_{ij}$ and the nominal link capacity $\kappa_{ij}$ are provided along with the network data.
	Considering the original trip demand $q$, we tested another two demand levels, 
	$1.5q$ (540,900 trips) and $2.0q$ (1,081,800 trips), to examine how the convergence processes are affected by the congestion level of the network. 
	
	To check the convergence of the NGEV equilibrium assignment, for the link flow and cost, we define the values at convergence $(\mathbold{X}^*, \mathbold{c}^*)$ and at the $m$-th iteration $(\mathbold{X}^{(m)}, \mathbold{c}^{(m)})$. We then define the relative differences $\eta_x \equiv \max_{ij \in \mathcal{L}} \frac{|X^{(m)}_{ij} - X^*_{ij}|}{X^*_{ij}}$ and $\eta_c \equiv \max_{ij \in \mathcal{L}} \frac{|c^{(m)}_{ij} - c^*_{ij}|}{c^*_{ij}}$ for the actual convergences of the problems \textbf{[NGEV-FD/P]} and \textbf{[NGEV-FD/D]}, respectively.\\
	
	
	\noindent\textbf{Primal Algorithms Convergence}
	
\begin{figure}[htb]
		\begin{center}
			\includegraphics[width=10cm]{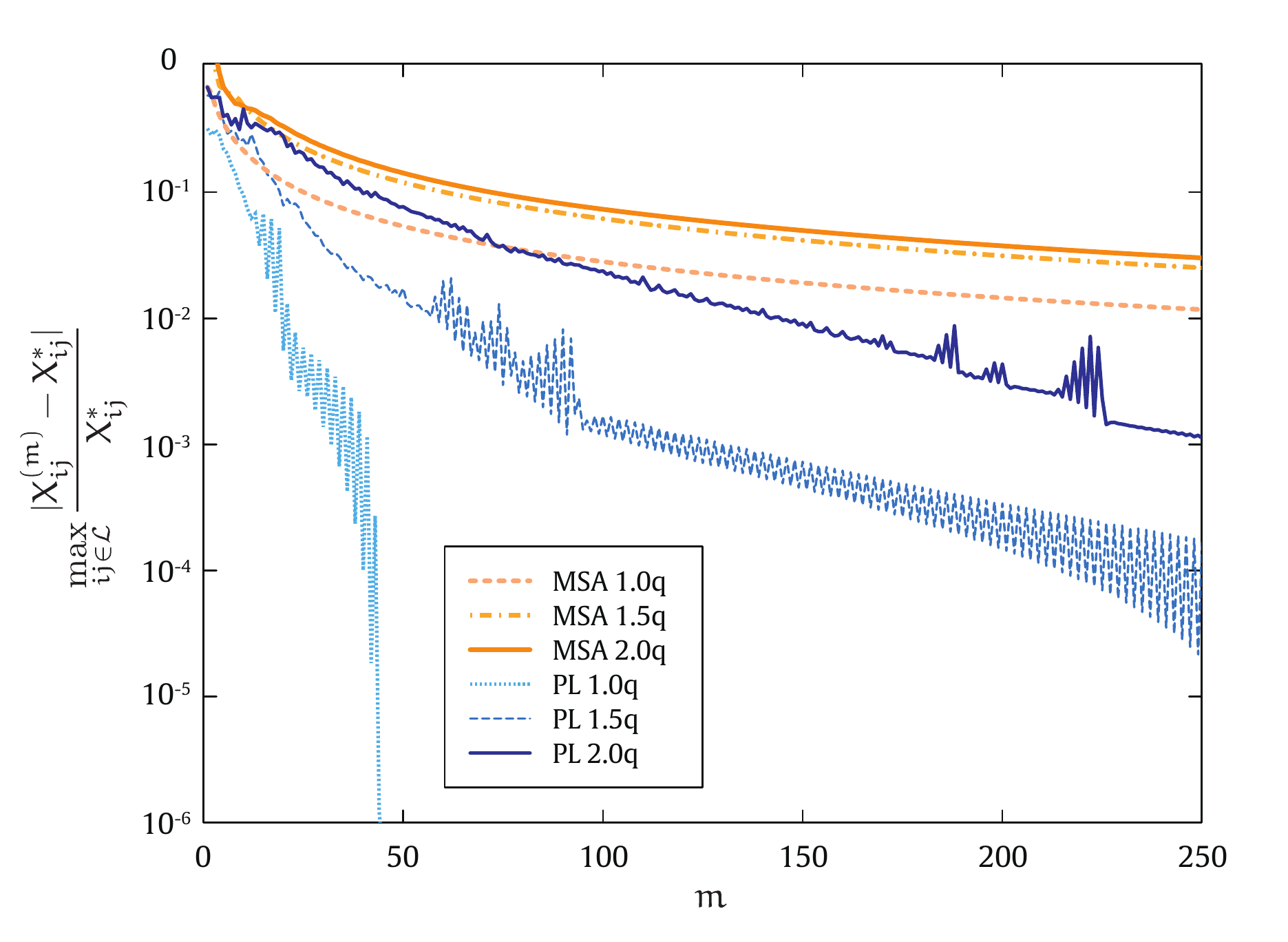}
			\caption{Convergence processes of MSA and PL with different demand levels.} 
			\label{fig:xdif} 
		\end{center}
	\end{figure}

	\noindent 
	Fig. \ref{fig:xdif} shows the convergence processes of the primal algorithms, MSA as a baseline, and PL, with the three different demand levels $\{q, 1.5q, 2.0q\}$. We observe that both algorithms properly move toward convergence with all demand levels. However, MSA is not practically usable because it is slow to converge, and $\eta_x$ remains above $10^{-2}$ even after 250 iterations. 
	In contrast, PL always converges more quickly and achieves considerably smaller values of $\eta_x$. However, the influence of the demand level on the convergence is not negligible. Whereas PL achieves $\eta_x = 10^{-6}$ before 50 iterations in the case with $q$, it does not achieve the same level of accuracy even after 250 iterations when the demand level is higher ($1.5q$ and $2.0q$ cases). \\
	
	\noindent \textbf{Dual Algorithms Convergence}
	
\begin{figure}[tbhp]
		\begin{center}
			\includegraphics[width=10cm]{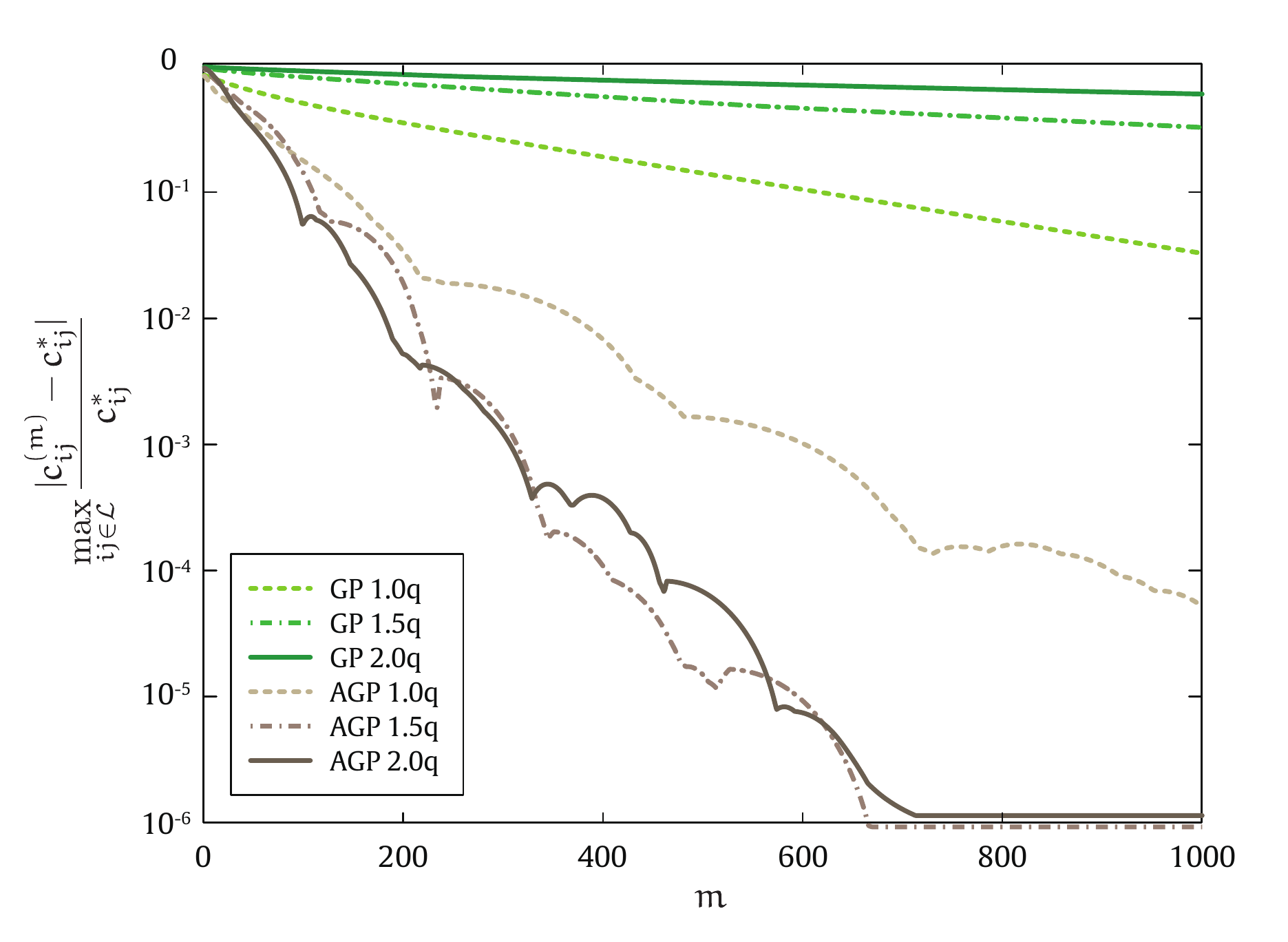}
			\caption{Convergence processes of GP and AGP with different demand levels.} 
			\label{fig:cdif} 
		\end{center}
	\end{figure}

	\noindent
	Fig. \ref{fig:cdif} shows the convergence processes of the dual algorithms, GP as a baseline, and AGP, with the three different demand levels $\{q, 1.5q, 2.0q\}$. For GP, we have to choose the step size $s$. After some trials, we set it in this experiment to $s = 10^{-5}$.
	For AGP, we set $k_{\min} = 50$ and apply the backtracking procedure with $\xi = 0.95$. The effectiveness of backtracking is proved though a sensitivity analysis in \ref{app:sensitivity}.
	In Fig. \ref{fig:cdif}, although both algorithms seem to move toward convergence, the improvements of GP are substantially slower. AGP converges quickly, and the convergence speed seems to be unaffected by the demand level. This result clearly shows the efficiency of the AGP method, i.e., the modification of the updating phase by a momentum significantly improves the computational performance.
	
	\subsection{Primal-Dual Algorithms Comparison}\label{sec:comparison}\noindent
	We then show a comparison between the algorithms for both primal and dual problems together. The same network data is used as the previous experiment, and the hyperparameters of the AGP method (with backtracking) are set to $k_{\min} = 50$ and $\xi = 0.25$.
	
\begin{figure}[tbhp]
		\begin{center}
			\includegraphics[width=10cm]{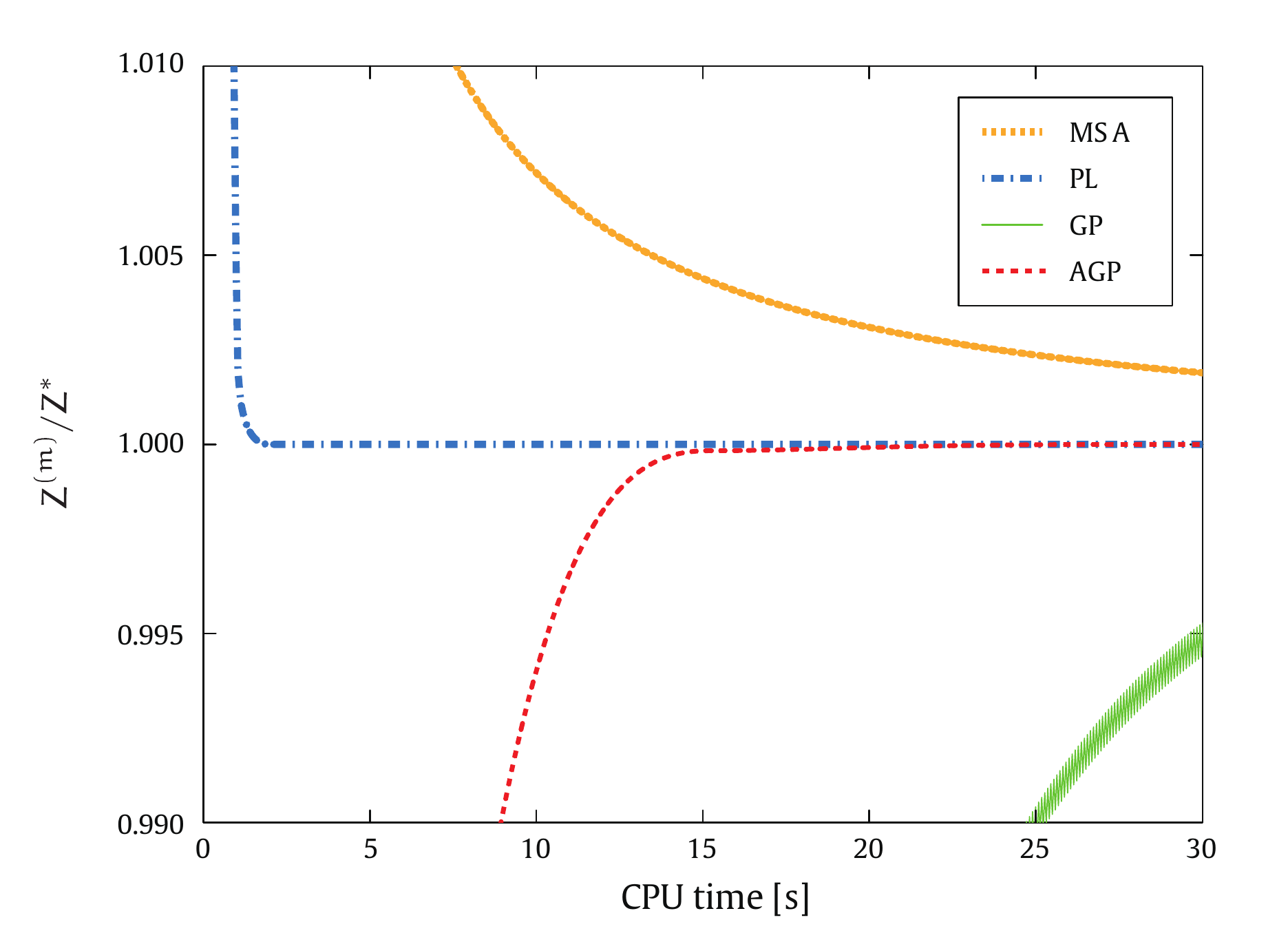}
			\caption{Comparison of the primal and dual algorithms with congestion level $q$.} 
			\label{fig:objective} 
		\end{center}
		
		\begin{center}
			\includegraphics[width=10cm]{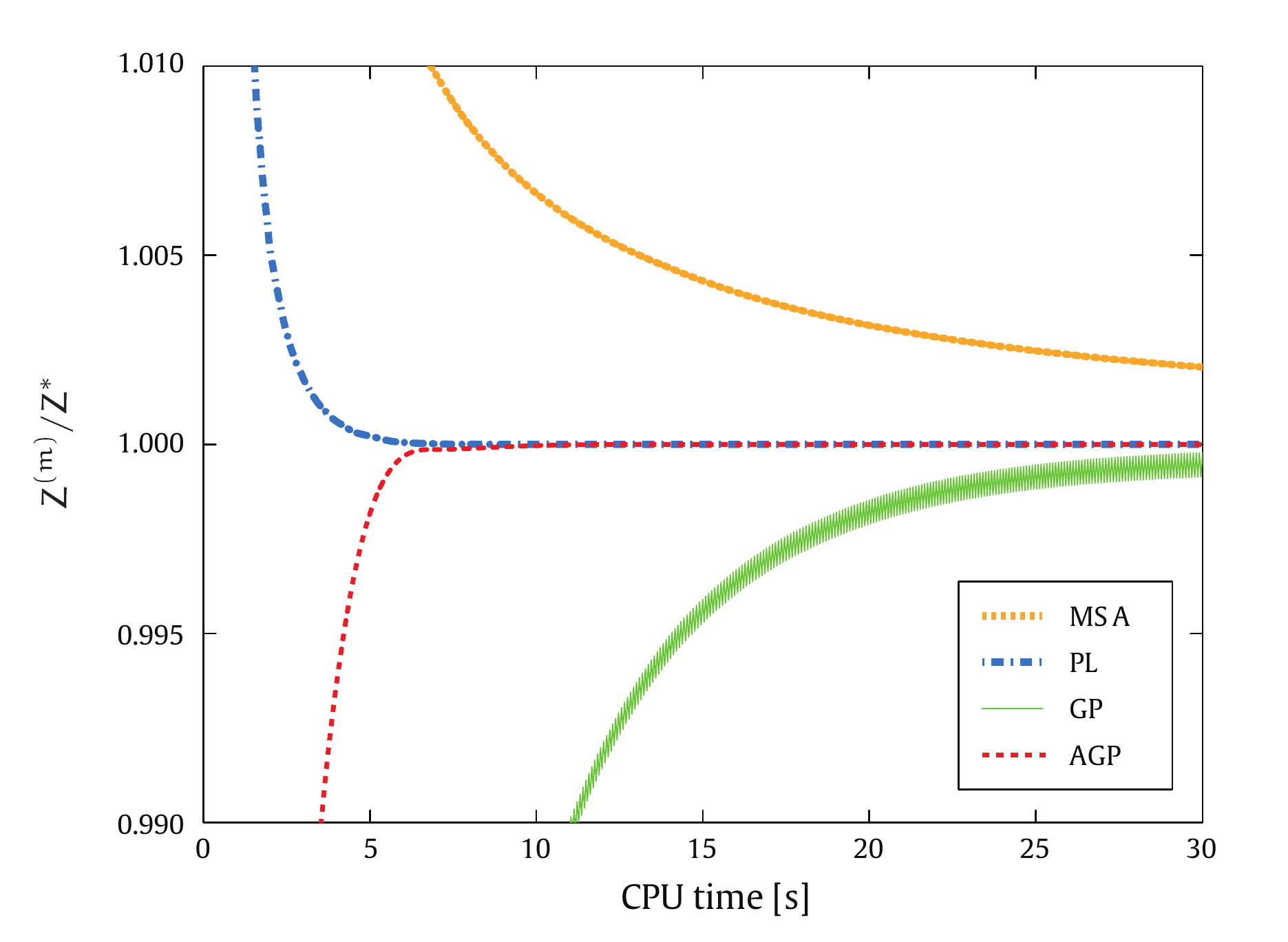}
			\caption{Comparison of the primal and dual algorithms with congestion level $1.5q$.} 
			\label{fig:objective_1.5q} 
		\end{center}
		
		\begin{center}
			\includegraphics[width=10cm]{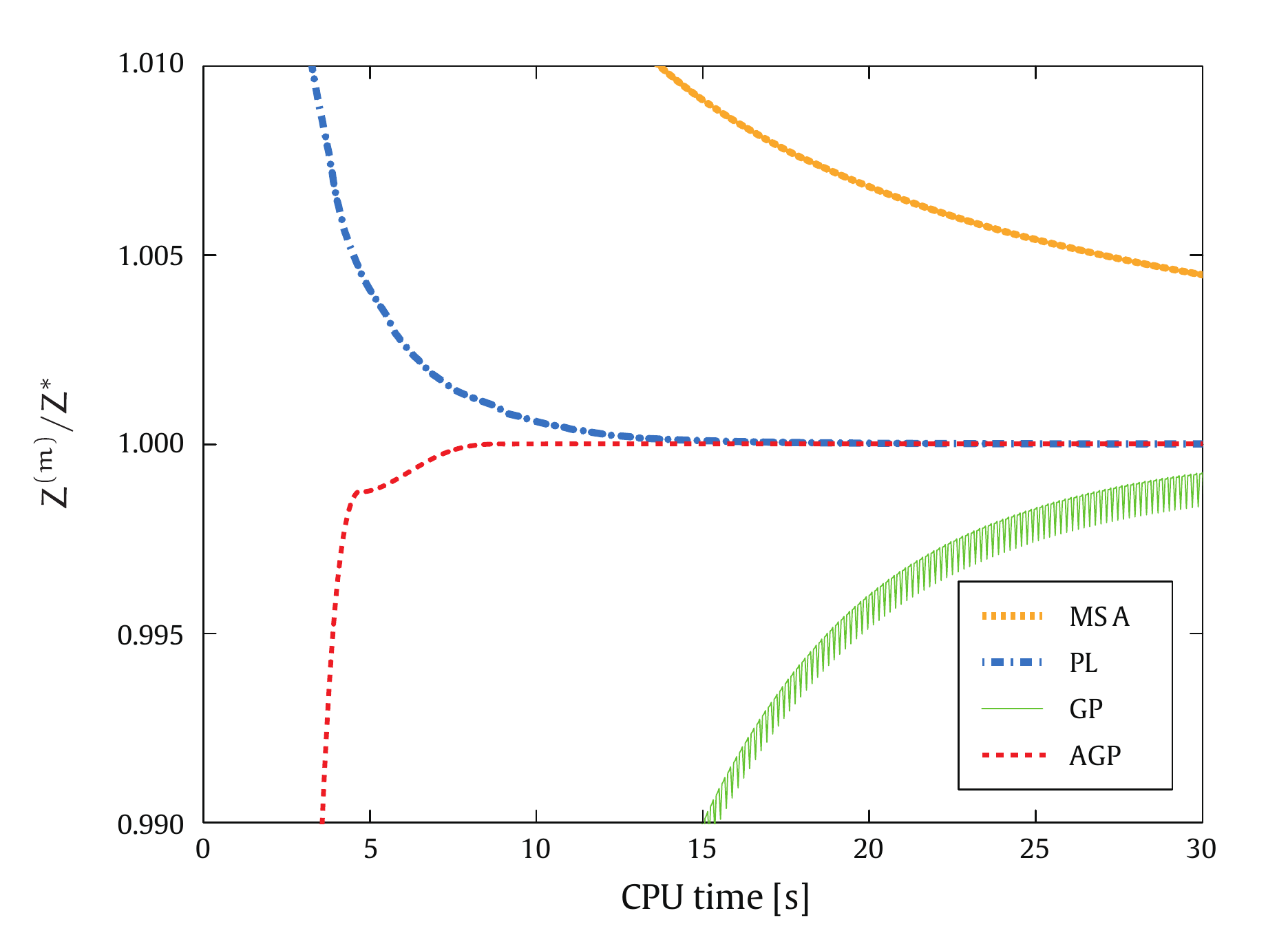}
			\caption{Comparison of the primal and dual algorithms with congestion level $2.0q$.} 
			\label{fig:objective_2q} 
		\end{center}
	\end{figure}

	For the sake of comparison we focus here on the objective value and its trajectory to the convergence for each algorithm because the optimal values of the objectives ($Z^*$) of the primal and dual problems are in theory consistent. 
	In fact, PL and AGP achieved the same objective values with $10^{-10}$ or a smaller order of the relative difference $\frac{|Z^{\rm FD *}_{\rm P} - Z^{\rm FD *}_{\rm D}|}{Z^{\rm FD *}_{\rm D}}$ for all the demand levels.
	The results are shown in Figs. \ref{fig:objective}-\ref{fig:objective_2q} respectively for each demand level of $\{q, 1.5q, 2.0q\}$, where the x-axis is set to the CPU time to compare the algorithms in terms of the computational efficiency. 
	
	With the original demand level $q$ (in Fig. \ref{fig:objective}), PL is the fastest to achieve the optimal objective value. AGP is the second..
	As the demand level increases, i.e., the network becomes more congested, the dual algorithms GP and AGP show their advantages. 
	In the case with $1.5q$ (in Fig. \ref{fig:objective_1.5q}), AGP is as efficient as PL, and GP slightly outperforms MSA. With $2.0q$ (in Fig. \ref{fig:objective_2q}), the greater performance of the dual algorithms is clearer. AGP is more efficient than the other algorithms. For all the demand levels, both PL and AGP show good convergence and clearly outperform the benchmarks MSA and GP.
	
	\subsection{Real-Size Application}\label{sec:realsize}
\begin{figure}[tb]
		\begin{center}
			\includegraphics[width=10cm]{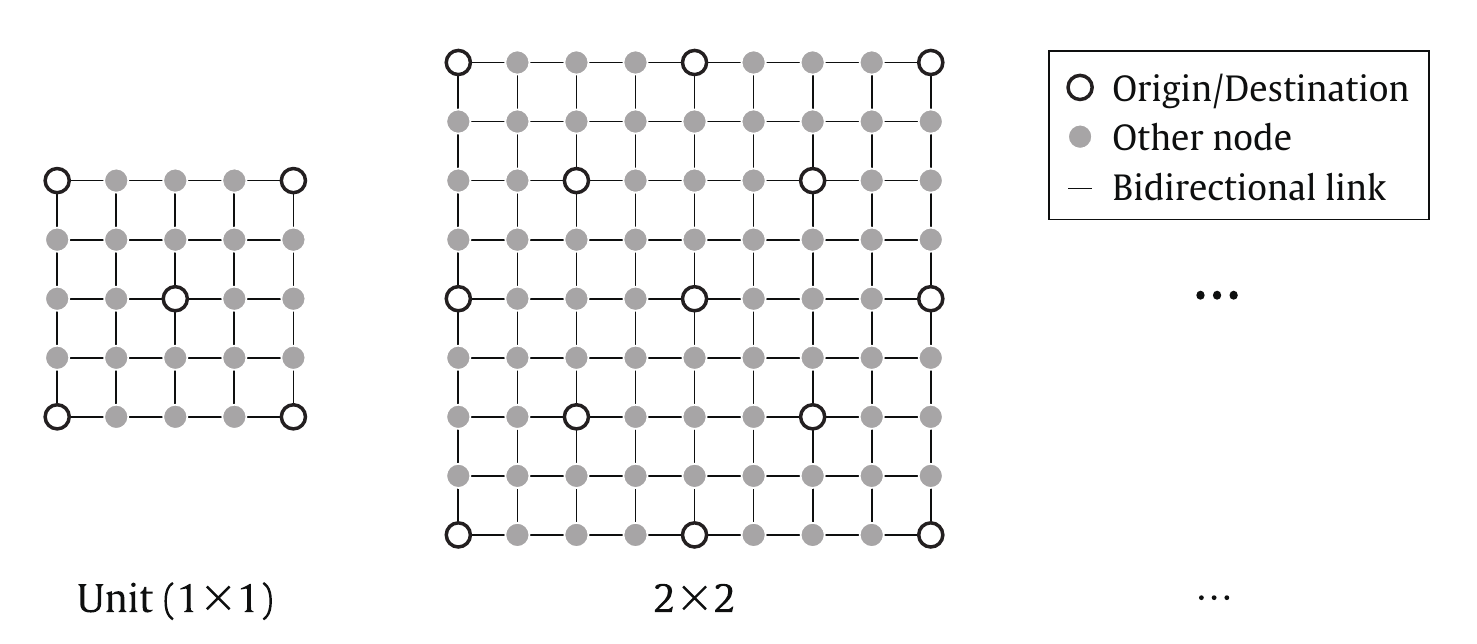}
			\caption{Bidirectional grid network, composed of $k \times k$ unit networks, where open circles indicate the origin/destination nodes. The $k \times k$ network has $|\mathcal{N}| = (4k + 1)^2$ nodes, $|\mathcal{L}| = 16k(4k+1)$ links, $|\mathcal{D}| = 1 + \sum_{l=1}^k 4k$ destinations, and $|\mathcal{W}| = |\mathcal{D}| \cdot (|\mathcal{D}|-1)$ OD pairs.}
			\label{fig:largenet} 
		\end{center}
	\end{figure}

\noindent
	We finally assess the computational efficiency of the proposed algorithms, PL and AGP, in real-size networks. We here do not consider MSA and GP as they clearly showed poor convergence even in a small network like the Sioux-Falls network.
	Bidirectional grid networks displayed in Fig.\ref{fig:largenet} are used for this experiment. We set for all $ij \in \mathcal{L}$: $\underline{c}_{ij} = 1$ and $\kappa_{ij} = 10000$. From each origin $o \in \mathcal{O}$, the generating flow proportional to its outdegree is assumed, i.e., $|\mathcal{F}(o)| \cdot q$ where $q$ is the reference flow. The OD flows are then defined by the following gravity model with $\nu = 0.1$:
	\begin{equation*}
	q^d_o = |\mathcal{F}(o)| \cdot q \cdot \frac{\exp(-\nu \underline{c}_{od})}{\sum_{d \in \mathcal{D}\backslash o} \exp(-\nu \underline{c}_{od})}.
	\end{equation*}
	The number of OD pairs is hence $|\mathcal{W}| = |\mathcal{D}| \cdot (|\mathcal{D}|-1)$. 
	We set the AGP hyperparameters at $k_{\min} = 50$ and $\xi = 0.25$, which are the same as the previous experiment.
	
	
\begin{figure}[tb]
		\begin{center}
			\includegraphics[width=10cm]{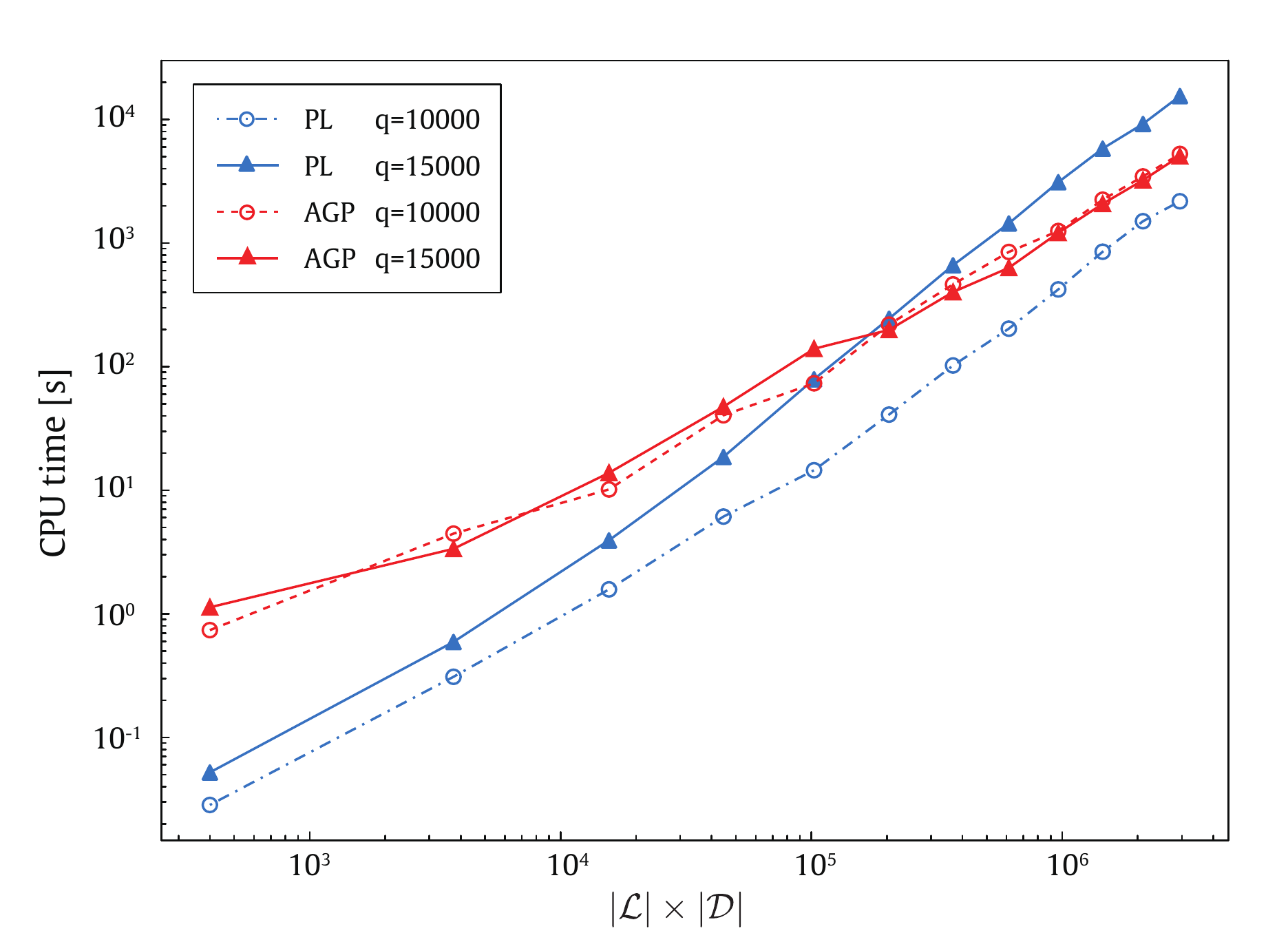}
			\caption{CPU times of NGEV equilibrium assignment in different sizes of the grid network.}
			\label{fig:gridtest} 
		\end{center}
	\end{figure}

	Fig.\ref{fig:gridtest} reports the CPU times that were required for the algorithms to converge. We defined the convergence with the criterion $\frac{|Z^{(m)} - Z^*|}{Z^*} \le 10^{-5}$, and tested in different sizes of the grid network with two different values of the reference flow $q = \{10000, 15000\}$. Fig.\ref{fig:gridtest} shows that the CPU time for each algorithm is approximately linear in the number of links multiplied by the number of destinations, $|\mathcal{L}| \times |\mathcal{D}|$. We again observe a clear difference by demand level in the results of PL. With $q=15000$ PL is clearly inefficient, and in the largest network ($k=12$), it requires over $10^4$ seconds to converge. In contrast, the effect of demand level on AGP is not significant. More importantly, AGP demonstrates its efficiency as the network size grows. The increase in CPU time of AGP is less rapid than that of PL, and AGP clearly outperforms PL in large-scale congested networks. 
	
	
	\subsection{Discussion}\noindent
	To summarize, the numerical examples showed clearly better convergence of PL and AGP, which we proposed in this study, as compared to the benchmarks MSA and GP. 
	Moreover, the primal-dual algorithm comparison clarified their complementary relationship. PL was more efficient than AGP in the case with the original demand level, but its performance declined with an increase in the demand level. In contrast, AGP was not affected by the demand level, and it outperformed PL when the network had higher demand levels. This comes from the difference in the unknown variables that the algorithms seek to find: link flow pattern $\mathbold x$ in the primal algorithm and link cost pattern $\mathbold c$ in the dual algorithm. 
	When the link flow approaches the link capacity, the link cost, as the return of the performance function, is sensitive to the change of flow. Therefore, in congested networks, the primal algorithms suffer from poor convergence, and the dual algorithms directly solving the link cost as the unknown show better convergence.
	
	Furthermore, the application to real-size networks showed the efficiency of AGP. Although PL is more efficient in less congested networks, the increase in computational effort of AGP by network size is less rapid than that of PL. This result again comes from the difference in the unknown variables. The solution space of AGP is $\mathbold{c} \in \mathbb{R}_+^{|\mathcal{L}|}$, whereas that of PL is $\mathbold{x} \in \mathbb{R}_+^{|\mathcal{D}| \times |\mathcal{L}|}$. That is, the number of unknowns of PL does not depend only on the number of links but also on the number of destinations. This difference may be crucial in the application to real networks because such networks often contain a high number of OD pairs as well as link/node sizes. Note that the computational efficiency of AGP is also affected by the number of destinations as the MTA loading procedure is destination-specific, but the same is true of PL.
	
	
	\section{Concluding Remarks}\label{sec:conclusion}\noindent
	This paper presented a framework of Markovian traffic equilibrium assignment, i.e., the NGEV equilibrium assignment. 
	The NGEV route choice model that captures the path correlation without explicit path enumeration was recently proposed by \cite{Hara2012NGEV, Hara2014NGEV} and \cite{Papola2013NGEV}. 
	However, its theoretical properties in the traffic assignment remained uninvestigated. The implementation of the NGEV route choice model in traffic equilibrium assignment requires the analysis on the solution properties and a convergent and efficient assignment algorithm, and this study provided the theoretical developments to this end.. 
	
	We first introduced and reformulated the NGEV route choice model to connect to the link-based traffic assignment formulation of \cite{Akamatsu1997Entropy}. Noticing that the NGEV assignment shares a common algebraic structure with the SP and MTA models \citep{Floyd1962, Warshall1962, Akamatsu1996MCA, Fosgerau2013RL}, we introduced a path algebra based on the work by \cite{Carre1979graphs}. The algebra provided a unified perspective to the MTA models including the NGEV assignment, and it was shown that any MTA model can be formulated and solved as the same system of equations under the path algebra.
	The compatibility of the NGEV assignment with MTA loading was shown through an illustrative example. We also showed that the NGEV assignment with MTA could alleviate the limitations of both the Dial-based NGEV assignment and the logit-based assignment. It is also important that the MTA fixes the path set and ensures the convergence to a consistent solution of the NGEV equilibrium assignment.
	
	We then provided the definition and formulations of the NGEV equilibrium assignment. The equivalent optimization problem was formulated by extending the link-based stochastic traffic equilibrium assignment model of \cite{Akamatsu1997Entropy}. Based on the formulation, we analyzed the uniqueness of the solution of the NGEV equilibrium assignment. Moreover, we presented the Lagrangian dual formulation that solves the NGEV equilibrium assignment in the space of link cost variable. The gradient of the dual objective function can be computed as the difference between the link demand function (i.e., flow-independent NGEV assignment with the current link cost pattern) and the inverse link cost function.
	
	We finally proposed solution algorithms for both the primal and dual formulations of the NGEV equilibrium assignment. The primal algorithm is based on the link-based PL method, which efficiently optimizes the step size and outperforms a naive iterative algorithm like MSA \citep{Akamatsu1997Entropy, Lee2010pl}. The dual algorithm is based on the accelerated gradient methods that have been developed recently in the machine learning field to efficiently solve large-scale optimization problems \citep{Nesterov1983, Beck2009fista}. This study was the first to investigate the use of the accelerated gradient methods in the traffic assignment field. In fact, the algorithm can be applied to logit-based and other traffic equilibrium assignment models. 
	The numerical experiments demonstrated the excellent convergence of the proposed algorithms. 
	Moreover, the complementary relationship between the primal and dual algorithms was observed. The PL method was more efficient than the AGP method in moderately congested networks, but the performance of PL declined with an increase in the demand level. In contrast, AGP was unaffected by the demand level, and it outperformed PL in congested and large-scale networks. 
	These results suggest that the proposed algorithms can be used in practice and selected depending on the network size and congestion level.
	
	To conclude, this paper provided a set of necessary theoretical developments and analyses for the NGEV equilibrium assignment: the generalized assignment framework unifying the MTA models under the same algebra, the equivalent optimization formulations, and the efficient solution algorithms. 
	Our contributions open up the applicability of NGEV-based models in the traffic assignment field. 
	The framework can be applied not only to vehicular networks, but also to transit networks based on frequency \citep[e.g.,][]{Lam1999stochastic, Ma2015NGEV} or with a timetable \citep[e.g.,][]{Nielsen2006transit, Nuzzolo2001doubly}, and to other various types of networks.
	Whereas this paper focuses on theoretical analyses, additional developments toward the practical application have remained for future research. 
	First, the development of the methods to estimate and validate the NGEV assignment models with real observed data is a major issue to be addressed in future research. This paper borrowed parameter values from the literature, but it is a remaining challenge to answer if the setting is transferable to any network.
	Although estimations of the MTA models with observed paths \citep[e.g.,][]{Fosgerau2013RL, Mai2015NRL, Mai2016NGEV} or with raw trajectory data \citep[e.g.,][]{Oyama2018pedestrian, Oijen2020wifi} have been presented, that of the NGEV equilibrium assignment additionally requires network-wide data such as link flow counts and OD matrices. The NGEV assignment comprises a large number of parameters, and technical developments of efficient estimation and data fusion are highly appreciated. 
	Second, the network reduction methods for the MTA models \citep[e.g.,][]{Oyama2019prism, Kazagli2020operational} can be integrated into the proposed framework. It makes the NGEV assignment more operational, which may be significant in application to large-scale networks.
	
	\section*{Acknowledgements}\noindent
	This study was supported by JSPS KAKENHI, Grant Numbers 18H01551 and 20K14899. 
	
	
	
	\appendix
	
	\section{List of notations} \label{app:notation}\noindent
	For the convenience of readers, below we list notations frequently used in this paper.\\
	
	\begin{longtable}[l]{@{}ll@{}} 
		\textbf{Sets} &  \\
		$\mathcal{N}$ & Set of nodes \\
		$\mathcal{L}$ & Set of links \\
		$\mathcal{O}$ & Set of origins $\subseteq {\mathcal N}$\\
		$\mathcal{D}$ & Set of destinations $\subseteq {\mathcal N}$\\
		$\mathcal{W}$ & Set of origin-destination pairs \\
		$\mathcal{F}(i)$ & Set of successor nodes of node $i$ \\
		$\mathcal{B}(i)$ & Set of predecessor nodes of node $i$ \\
		$\mathbb{R}_{\epsilon}$ & Set of elements for the path algebra $\mathbb{R}_{\text{path}}$ \\
		$\mathcal{M}_{n}(\mathbb{R}_{\epsilon})$ & Set of $n \times n$ matrices whose entries belong to $\mathbb{R}_{\epsilon}$ \\
		\\
		\textbf{Parameters} &  \\
		$\underline{c}_{ij}$ & Free-flow cost of link $ij \in \mathcal{L}$\\
		$\kappa_{ij}$ & Capacity of link $ij \in \mathcal{L}$\\
		$q^d_i$ & Given demand flow between node $i \in \mathcal{N}$ and destination $d \in \mathcal{D}$ \\
		$\tilde{q}^d_i$ & Modified demand flow between node $i \in \mathcal{N}$ and destination $d \in \mathcal{D}$ \\
		$\theta^d_i$ & NGEV scale parameter for node $i \in \mathcal{N}$ and destination $d \in \mathcal{D}$ \\
		$\alpha^d_{ji}$ & NGEV allocation parameter of node $j \in \mathcal{N}$ to $i \in \mathcal{B}(j)$ by destination $d \in \mathcal{D}$ \\
		$s$ & Step size for the dual algorithms\\
		$k_{\min}$ & Minimum number of iterations for adaptive restart for AGP\\
		$\xi$ & Rate of updating step size during backtracking for AGP\\
		\\
		\textbf{Variables} &  \\
		$\mu^d_{i}$ & Expected minimum cost from node $i \in \mathcal{N}$ to destination $d \in \mathcal{D}$\\
		$p^d_{ij|i}$ & Link choice probability defined for node pair $i,j \in \mathcal{N}$ and destination $d \in \mathcal{D}$\\
		$x^d_{ij}$ & Flow on link $ij \in \mathcal{L}$ directed toward destination $d \in \mathcal{D}$ \\
		$z^d_{i}$ & Inflow to node $i \in \mathcal{N}$ directed toward destination $d \in \mathcal{D}$\\
		$X_{ij}$ & Flow on link $ij \in \mathcal{L}$ \\
		$c_{ij}(X_{ij})$ & Cost function of link $ij \in \mathcal{L}$ \\
		$c_{ij}^{-1}(c_{ij})$ & Inverse of the link cost function of link $ij \in \mathcal{L}$ \\
		$w_{ij}$ & Weight of link $ij \in \mathcal{L}$ for the path algebra $\mathbb{R}_{\text{path}}$ \\
		$H^d_i(\mathbold{x}^d)$ & Entropy function for node $i \in \mathcal{N}$ and destination $d \in \mathcal{D}$\\
		$\eta_x$ & Relative flow difference index for convergence \\ 
		$\eta_c$ & Relative cost difference index for convergence \\ 
	\end{longtable}
	
	\section{Proofs} \label{sec:proof}\noindent
	\subsection{Proof of Proposition 1}\noindent
	The total cost function $\mathbold{c} \cdot \mathbold{x}^d$ is convex, and the entropy function (\ref{eq:entropy}) is strictly concave. Because the objective function (\ref{eq:ZNGEV}) is the sum of a convex function and a strictly convex function, it is strictly convex\footnote{We omit the detailed discussion, but the strict convexity can be easily proved by examining the elements of Hessian matrix. \cite{Akamatsu1997Entropy} provided a detailed proof for the case of logit-based assignment, and our proof is its straightforward extension.}.
	Moreover, the problem \textbf{[NGEV]} has only linear constraints and a non-negative condition. Hence, the feasible region of the problem is a closed, nonempty convex set.
	The problem \textbf{[NGEV]} is thus a convex programming problem, and the Karush-Kuhn-Tucker (KKT) condition is the necessary and sufficient condition for the optimality.
	
	The Lagrangian of \textbf{[NGEV]} is given by
	\begin{equation}
	L(\mathbold{x}, \bm{\mu})
	\equiv
	Z(\mathbold{x}) + \sum_{d \in \mathcal{D}} \bm{\mu}^d \cdot (\tilde{\mathbold{q}}^d - \mathbf{A}\mathbold{x}^d),
	\label{eq:Lngev}
	\end{equation}
	where $\bm{\mu}$ is the Lagrangian multiplier.
	In the stochastic assignment based on the RUMs, every alternative (link) is always chosen with a strictly positive (non-zero) probability. Therefore, the optimal flow must be positive, i.e., $\mathbold{x}^* > 0$, and the first-order condition is given by:
	\begin{equation}
	\mathbold{x}^* 
	\equiv
	\arg\min_{\mathbold{x} \ge 0} L(\mathbold{x}, \bm{\mu})
	\Leftrightarrow
	\mathbold{0} \le \mathbold{x}^* \perp \nabla_{\mathbold x} L(\mathbold{x}^*)
	\Leftrightarrow
	\nabla_{\mathbold x} L(\mathbold{x}^*) = 0,
	\end{equation}
	and for all $ij \in \mathcal{L}$ and $d \in \mathcal{D}$, the derivative of Lagrangian with respect to the link flow is:
	\begin{equation}
	\frac{\partial L}{\partial x^d_{ij}} 
	= \frac{1}{\theta_i} \left[\ln \frac{x_{ij}^{d}}{\alpha^d_{ji}z_i^{d}} + \theta_i (c_{ij} + \mu^d_j - \mu^d_i) \right].
	\end{equation}
	Considering (\ref{eq:zi}), we finally obtain
	\begin{equation}
	p^{d*}_{ij|i} \equiv \frac{x_{ij}^{d*}}{\sum_{j \in \mathcal{F}(i)} x^{d*}_{ij}}  = \alpha^d_{ji} e^{-\theta^d_i (c_{ij} + \mu^d_j - \mu^d_i)},
	\label{eq:pstar}
	\end{equation}
	and
	\begin{equation}
	1 = \sum_{j \in \mathcal{F}(i)} \alpha^d_{ji} e^{-\theta^d_i (c_{ij} + \mu^d_j - \mu^d_i)}.
	\end{equation}
	These results are consistent with the NGEV route choice model defined by (\ref{eq:pij}) and (\ref{eq:mu}), and the Lagrangian multiplier $\bm{\mu}$ corresponds to the expected minimum cost.
	\hfill $\square$
	
	\subsection{Proof of Corollary 1}\noindent
	From the discussion in the proof of \textbf{Proposition \ref{prop:DTA}}, the problem \textbf{[NGEV]} is a convex programming problem that has a strictly convex objective function. Thus, the globally optimal solution can be uniquely determined if a solution exists. 
	\hfill $\square$
	
	\subsection{Proof of Proposition 2}\label{app:proof_prp2} \noindent
	The only difference between \textbf{[NGEV]} and \textbf{[NGEV-FD/P]} is the first term of their objective functions. Based on \textbf{Assumption \ref{assumption:cX}} regarding the link cost function, we know that the first term of the objective function (\ref{eq:NGEVSUE}), i.e., (\ref{eq:Cint}), is strictly convex. Given the proof of \textbf{Proposition \ref{prop:DTA}}, this is sufficient to prove that \textbf{[NGEV-FD/P]} is a convex programming problem. Then, the KKT condition is the necessary and sufficient condition for the optimality.
	
	The Lagrangian of \textbf{[NGEV-FD/P]} is given by
	\begin{equation}
	L^{\rm FD}_{\rm P}(\mathbold{x}, \bm{\mu})
	\equiv
	Z^{\rm FD}_{\rm P}(\mathbold{x}) + \sum_{d \in \mathcal{D}} \bm{\mu}^d \cdot (\tilde{\mathbold{q}}^d - \mathbf{A}\mathbold{x}^d),
	\label{eq:LSUE}
	\end{equation}
	and $\nabla_{\mathbold x} Z^{\rm FD}_{\rm P} = (c_{ij}(X_{ij}))_{ij \in \mathcal{L}}$.
	In the same way as \textbf{[NGEV]}, the first-order condition provides the NGEV equilibrium assignment conditions (\ref{eq:pij})-(\ref{eq:edgeflow}), (\ref{eq:linkflow_vector}) and (\ref{eq:dcdX}).
	\hfill $\square$
	
	\subsection{Proof of Corollary 2}\noindent
	From the discussion in the proof of \textbf{Proposition \ref{prop:SUE}}, the problem \textbf{[NGEV-FD/P]} is a convex programming problem that has a strictly convex objective function. Thus, the globally optimal solution can be uniquely determined if a solution exists.
	\hfill $\square$
	
	\subsection{Proof of Proposition 3}\noindent
	By using (\ref{eq:pstar}), for all $i \in \mathcal{N}$ and $d \in \mathcal{D}$, the optimal entropy is obtained by
	\begin{align}
	H^d_i(\mathbold{x}^{d*})
	&= -\sum_{j \in \mathcal{F}(i)} x_{ij}^{d*} \ln \frac{p_{ij|i}^{d*}}{\alpha_{ji}} \nonumber\\
	&= \theta^d_i \left[\sum_{j \in \mathcal{F}(i)} c_{ij} x^{d*}_{ij} - \sum_{j \in \mathcal{F}(i)} (\mu^d_i - \mu^d_j) x^{d*}_{ij} \right],
	\label{eq:Hxstar}
	\end{align}
	and therefore, 
	
	\begin{equation}
	\hat{\bm{\theta}} \cdot \mathbold{H}^d(\mathbold{x}^{d*}) 
	= \mathbold{c} \cdot \mathbold{x}^{d*} - (\mathbf{A}^{\top}\bm{\mu}^d) \cdot \mathbold{x}^{d*} 
	= \mathbold{c} \cdot \mathbold{x}^{d*} - \bm{\mu}^d \cdot (\mathbf{A}\mathbold{x}^{d*}).
	\label{eq:Hvector}
	\end{equation}
	Substituting (\ref{eq:Hvector}) into the Lagrangian (\ref{eq:LSUE}) of \textbf{[NGEV-FD/P]} yields
	\begin{equation}
	L^{\rm FD}_{\rm P}(\mathbold{x}^*, \bm{\mu})
	= C(\mathbold{X}^*) - \mathbold{c} \cdot \mathbold{X}^* + \sum_{d \in \mathcal{D}}\bm{\mu}^d \cdot \tilde{\mathbold{q}}^d.
	\end{equation}
	From (\ref{eq:cstar}), with the optimal flow $\mathbold{X}^*$ the following holds:
	\begin{equation}
	\mathbold{c} \cdot \mathbold X^* -
	C(\mathbold{X}^*) = C^*(\mathbold{c}).
	\end{equation}
	Consequently, we obtain the dual problem of \textbf{[NGEV-FD/P]}, namely \textbf{[NGEV-FD/D]}:
	\begin{equation}
	\max_{\bm{\mu}^d \in \mathcal{K}_d} \min_{\mathbold{x} \ge 0} L^{\rm FD}_{\rm P}(\mathbold{x}, \bm{\mu})
	= \max_{\bm{\mu}^d \in \mathcal{K}_d} L^{\rm FD}_{\rm P}(\mathbold{x}^*, \bm{\mu})
	= \max_{\mathbold{c}\ge\underline{\mathbold{c}}} 
	[
	- C^*(\mathbold{c})
	+ \sum_{d \in \mathcal{D}}\bm{\mu}^d(\mathbold c) \cdot \tilde{\mathbold{q}}^d
	].
	\end{equation}
	\hfill $\square$
	
	\subsection{Proof of Lemma 1}\noindent
	The derivative of $Z^{\rm FD}_{\rm D}$ is
	\begin{equation}
	\nabla_{\mathbold{c}} Z^{\rm FD}_{\rm D}(\mathbold{c}) 
	= \sum_{d \in \mathcal{D}} \frac{\partial \bm{\mu}^d(\mathbold c)}{\partial \mathbold{c}} \cdot \tilde{\mathbold{q}}^d 
	- \mathbold{c}^{-1}\mathbold{(c)}.
	\label{eq:grad_Z}
	\end{equation}
	By taking the logarithm of (\ref{eq:mu}), we have
	\begin{equation}
	\mu^d_i = -\frac{1}{\theta^d_i} \ln \sum_{j \in \mathcal{F}(i)} \alpha^d_{ji} e^{- \theta^d_i (c_{ij} + \mu^d_j)},
	\end{equation}
	whose derivative with respect to $c_{kl}$ is
	\begin{equation}
	\frac{\partial \mu^d_i}{\partial c_{kl}}
	= \sum_{j \in \mathcal{F}(i)} p^d_{ij|i} (
	\frac{\partial c_{ij}}{\partial c_{kl}}
	+ \frac{\partial \mu^d_j}{\partial c_{kl}}
	)
	= \delta_{ij,kl} + \sum_{j \in \mathcal{F}(i)} p^d_{ij|i} 
	\frac{\partial \mu^d_j}{\partial c_{kl}},
	\label{eq:grad_mu}
	\end{equation}
	as $\frac{\partial \mathbold c}{\partial \mathbold c}  = \mathbf I$, which is an identity matrix. Hence, with $\frac{\partial \mu^d_d}{\partial c_{kl}} = 0, \forall kl \in \mathcal{L}$,
	\begin{align}
	\sum_{i \in \mathcal{N}} \frac{\partial \mu^d_i}{\partial c_{kl}} q_{id}
	&= 
	\sum_{i \in \mathcal{N}} q_{id} \sum_{j'_1 \in \mathcal{F}(i)} p^d_{ij'_1|i} \sum_{j'_2 \in \mathcal{F}(j'_1)} p^d_{j'_1j'_2|j'_1} \cdots p^d_{j'_n k|j'_n} p^d_{kl|k} \nonumber\\
	&= \left( \sum_{i \in \mathcal{N}} q_{id} 
	\sum_{r \in \mathcal{R}_{ik}} p(r)
	\right) p^d_{kl|k} \nonumber\\
	&= z^d_k p^d_{kl|k} = x^d_{kl},
	\end{align}
	where $\mathcal{R}_{ik}$ is the set of all feasible routes between $i$ and $k$, and by substituting this into (\ref{eq:grad_Z}), the lemma is proved.
	\hfill $\square$
	
	\section{Computational time of network loading} \label{app:loadingtime}\noindent
	We report an experiment to compare the computational times of network loading algorithms. For the test, we used four different sizes of grid networks in Fig. \ref{fig:largenet} and set the same NGEV parameters as Model 3 in Section \ref{sec:loading_example}. Table \ref{tab:loading_time} reports the results. First and most importantly, the required CPU times of the NGEV assignment are in the same level as those of the logit assignment. As shown in Section \ref{sec:algebra}, both the logit and NGEV assignment can be performed by a link-based and many-to-one assignment procedure, meaning that the computation of expected minimum costs and link flows are solved only as many times as the number of destinations \citep[e.g.,][]{Mai2015NRL, Mai2016NGEV}. Note that, whereas Dial requires little additional effort to perform the NGEV assignment, there is a non-negligible difference between the results of logit- and NGEV-MTA. In logit-MTA, Eq.(\ref{eq:V}) reduces to a system of linear equations, which one can efficiently solve \citep{Akamatsu1996MCA, Fosgerau2013RL}. In NGEV-MTA, it becomes a nonlinear system, and Eq.(\ref{eq:V}) needs to be solved by the value iteration algorithm \citep[e.g.,][]{Mai2015NRL}. 
	Nevertheless, NGEV-MTA is still as efficient as Dial's algorithm, as shown in Table \ref{tab:loading_time}. Also, the logit and NGEV assignment models are incomparably faster than the probit model: at $k=8$, even a single draw for the probit (mean $= 1.104$s) is slower than logit-MTA ($0.910$s). Moreover, the probit assignment with a limited number of draws involves a significant approximation error (see \ref{app:probit}).
	
	\begin{table}[H]
		\centering 
		\footnotesize
		\caption{CPU times in seconds for network loading in different sizes of networks (the reported values are averages over 100 runs). The values in parentheses are ones for a single destination.}
		\label{tab:loading_time}
		\begin{tabular}{|c|c|c|c||c|c|c|c|c|}\hline
			$k$ & $|\mathcal{L}|$ & $|\mathcal{W}|$ & $|\mathcal{D}|$ & Logit-Dial & Logit-MTA & NGEV-Dial & NGEV-MTA & Probit ($100$ draws) \\
			\hline
			\multirow{2}{*}{1}             &  \multirow{2}{*}{80}       & \multirow{2}{*}{20}     & \multirow{2}{*}{5}        & 0.003     & 0.007    & 0.004    & 0.006 & 0.115 \\
			&         &      &         & (0.0006)     & (0.0014)    & (0.0008)    & (0.0012) & (0.0230) \\ \hline
			\multirow{2}{*}{2}             &  \multirow{2}{*}{288}      & \multirow{2}{*}{156}    & \multirow{2}{*}{13}       & 0.018     & 0.026    & 0.020    & 0.037 & 0.471 \\
			&         &      &         & (0.0014)     & (0.0020)    & (0.0015)    & (0.0028) & (0.0362) \\ \hline
			\multirow{2}{*}{4}            &  \multirow{2}{*}{1088}     & \multirow{2}{*}{1640}   & \multirow{2}{*}{41}       & 0.167     & 0.104    & 0.189    & 0.189 & 5.302 \\
			&         &      &         & (0.0041)     & (0.0025)    & (0.0046)    & (0.0046) & (0.1293) \\ \hline
			\multirow{2}{*}{8}            &  \multirow{2}{*}{4224}     & \multirow{2}{*}{20880}  & \multirow{2}{*}{145}      & 2.196     & 0.910    & 2.570    & 3.036 & 110.438 \\
			&         &      &         & (0.0151)     & (0.0063)    & (0.0177)    & (0.0209) & (0.7616) \\
			\hline
		\end{tabular}
	\end{table}
	
	\section{Flow patterns of logit and NGEV equilibrium assignment} \label{app:sue}
\begin{figure}[htbp]
		\begin{center}
			\includegraphics[width=15cm]{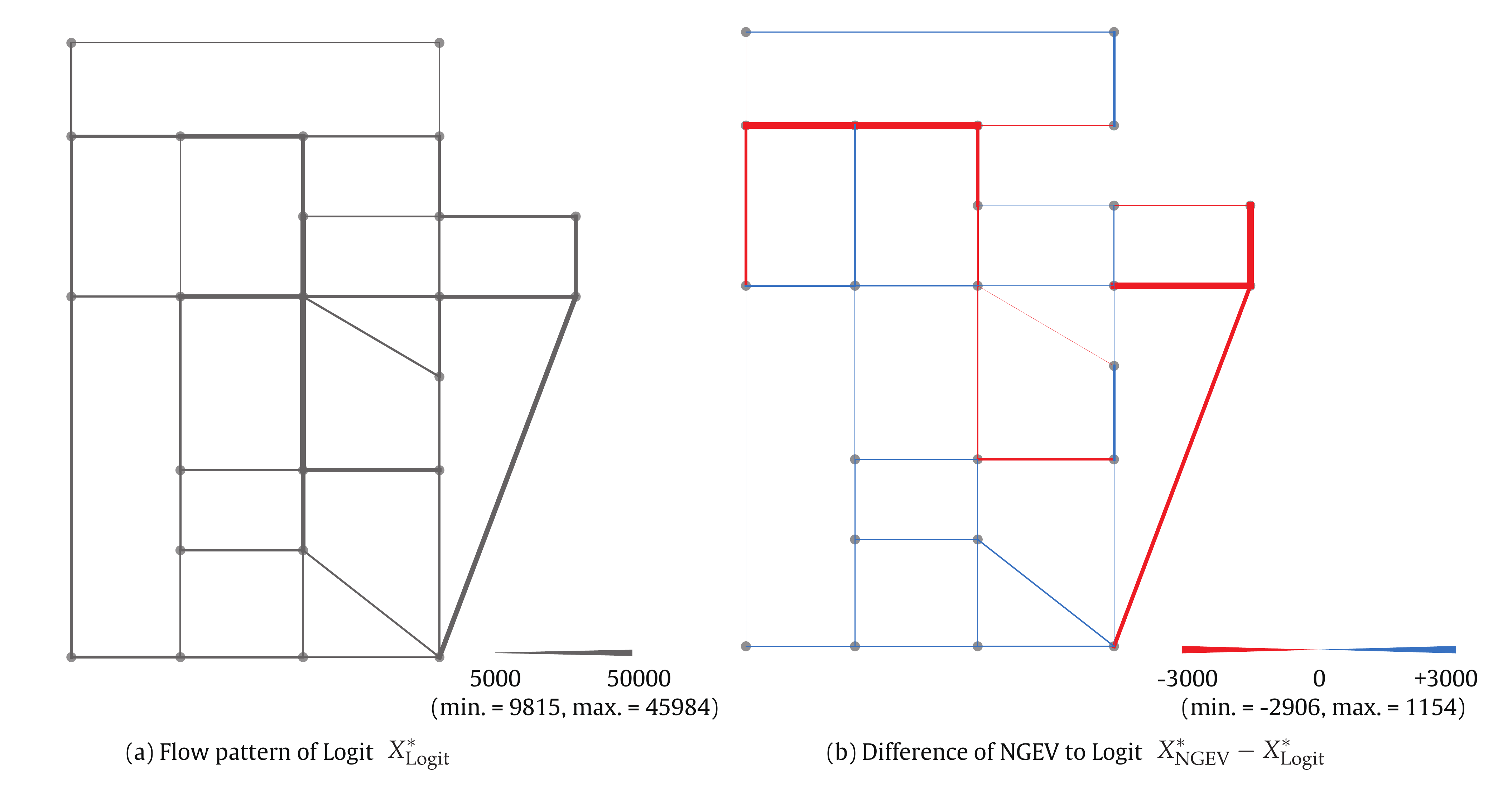}
			\caption{Equilibrated flow patterns. (a) Logit equilibrium assignment, and (b) the difference of NGEV equilibrium assignment to logit.} 
			\label{fig:SUE} 
		\end{center}
	\end{figure}

\noindent
	Fig.\ref{fig:SUE} shows the comparison of the results of logit and NGEV equilibrium assignment. The left panel (a) maps the equilibrated flow pattern $\mathbold{X}^*_{\text{Logit}}$ of logit ($\theta = 1.0$), and the right panel (b) displays the difference $\mathbold{X}^*_{\text{NGEV}} - \mathbold{X}^*_{\text{Logit}}$ of NGEV (Model 3 in Section \ref{sec:loading_example}) to logit. The NGEV model assigns fewer flows on links, displayed in red in (b), on which the logit model assigns a large volume of flows. This result shows that the NGEV model successfully captures the path correlation and reduces the utilities of overlapping paths. Note that we solved both the logit and NGEV equilibrium assignment models by the PL method with the convergence threshold of $0.1$\% error, and their runtimes were $3.85$s ($53$ iterations) and $2.87$s ($38$ iterations), respectively.
	
	With the equilibrated cost patterns, moreover, we analyze paths generated by logit and NGEV models. For each model, we sampled $10000$ paths departing from the top-left node and terminating at the bottom-right node. The logit model generated 11 unique paths, whereas the NGEV model generated $21$ ones. We also counted the number of generated cyclic paths that pass some nodes twice or more, and the logit and NGEV models respectively generated $156$ ($1.56$\%) and $115$ ($1.15$\%) cyclic paths out of $10000$. All the $115$ cyclic paths of the NGEV model contain only a single loop. The NGEV model increases the variety of paths and alleviates the generation of cyclic paths, although it may still generate cyclic paths. 
	
	\section{Incorporating Backtracking Procedure} \label{app:sensitivity}
\begin{figure}[htbp]
		\begin{center}
			\includegraphics[width=10cm]{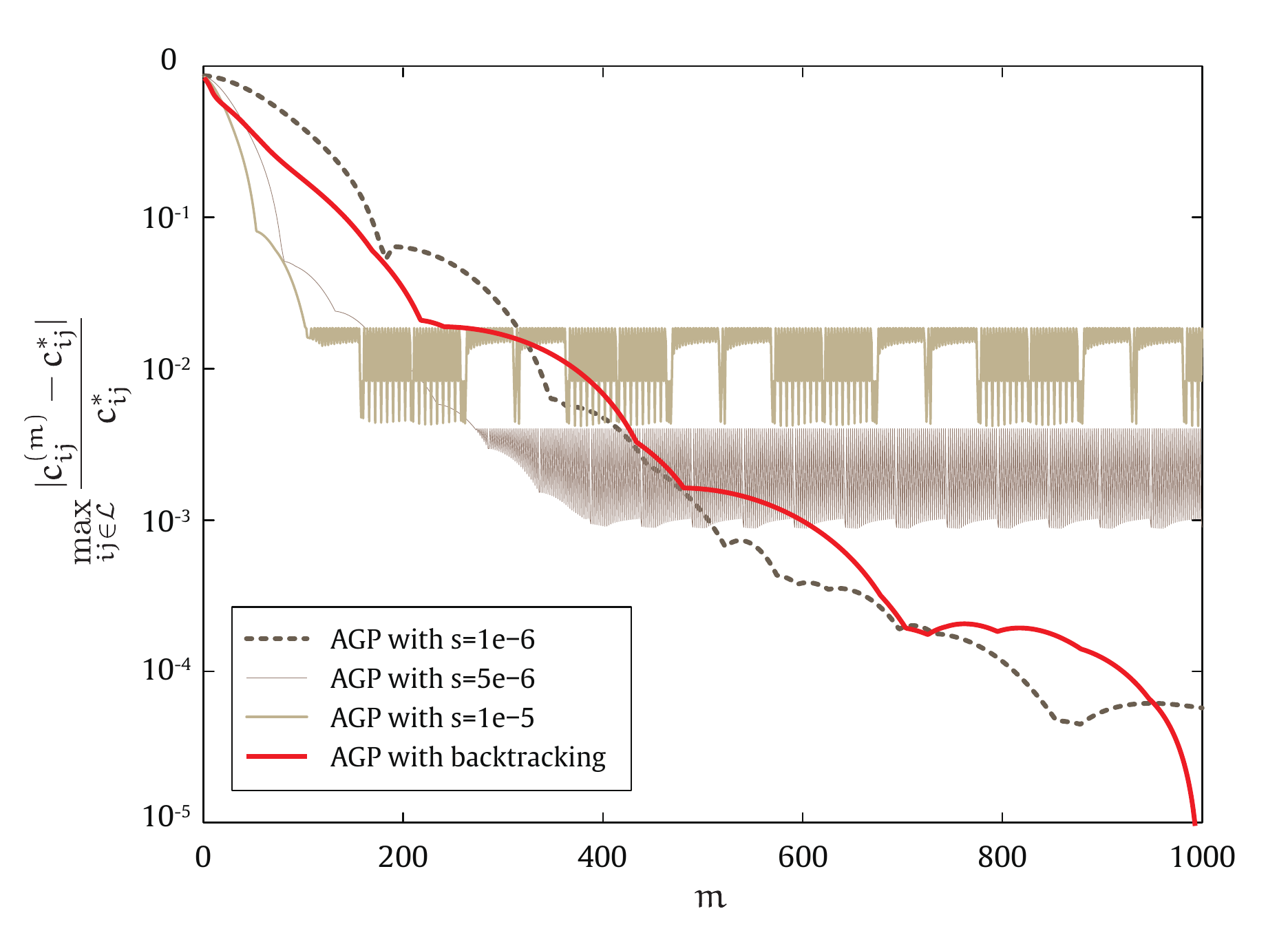}
			\caption{Convergence processes of AGP with different step sizes and with backtracking procedure.} 
			\label{fig:backtracking} 
		\end{center}
	\end{figure}

\noindent
	The AGP method relies on an arbitrarily chosen step size $s$, which influences the algorithm's performance. Incorporating the backtracking procedure into AGP provides a solution to this limitation. 
	In Fig. \ref{fig:backtracking}, we compare different step sizes for AGP and a case with backtracking, where the demand level is $q$ in the Sioux Falls network. 
	
	When we set $s$ to a relatively large value ($s=10^{-5}$ or $5\cdot 10^{-6}$), the solution improves quickly at the beginning, but later fluctuates, and $\eta_c$ does not become smaller than a certain level. A large step size may try updating the link cost outside of the feasible region, i.e., below the free-flow cost, which is not allowed and modified as (\ref{eq:projection}). This bounding causes the fluctuation observed in Fig. \ref{fig:backtracking}, also for GP in the experiment of Section \ref{sec:comparison}.
	A smaller step size ($s=10^{-6}$) works better in terms of stable convergence. The solution improvement is not as fast as the cases with larger step sizes but seems to be smooth until the convergence. 
	
	However, in any case, the step size $s$ is still arbitrary, and we will not know if there exist better values unless we conduct a number of trials. The backtracking procedure (\ref{eq:z_BT})-(\ref{eq:Q_BT}) solves this problem and achieves a fast and smooth solution improvement because it continues updating the step size $s$ during the iterations. As shown in Fig. \ref{fig:backtracking}, AGP with backtracking improves the solution more quickly than the case with a small step size ($s=10^{-6}$) at the beginning, and it then converges smoothly, achieving a higher quality solution (with $\eta_c < 10^{-5}$) than the other cases.
	
\begin{figure}[htbp]
		\begin{center}
			\includegraphics[width=10cm]{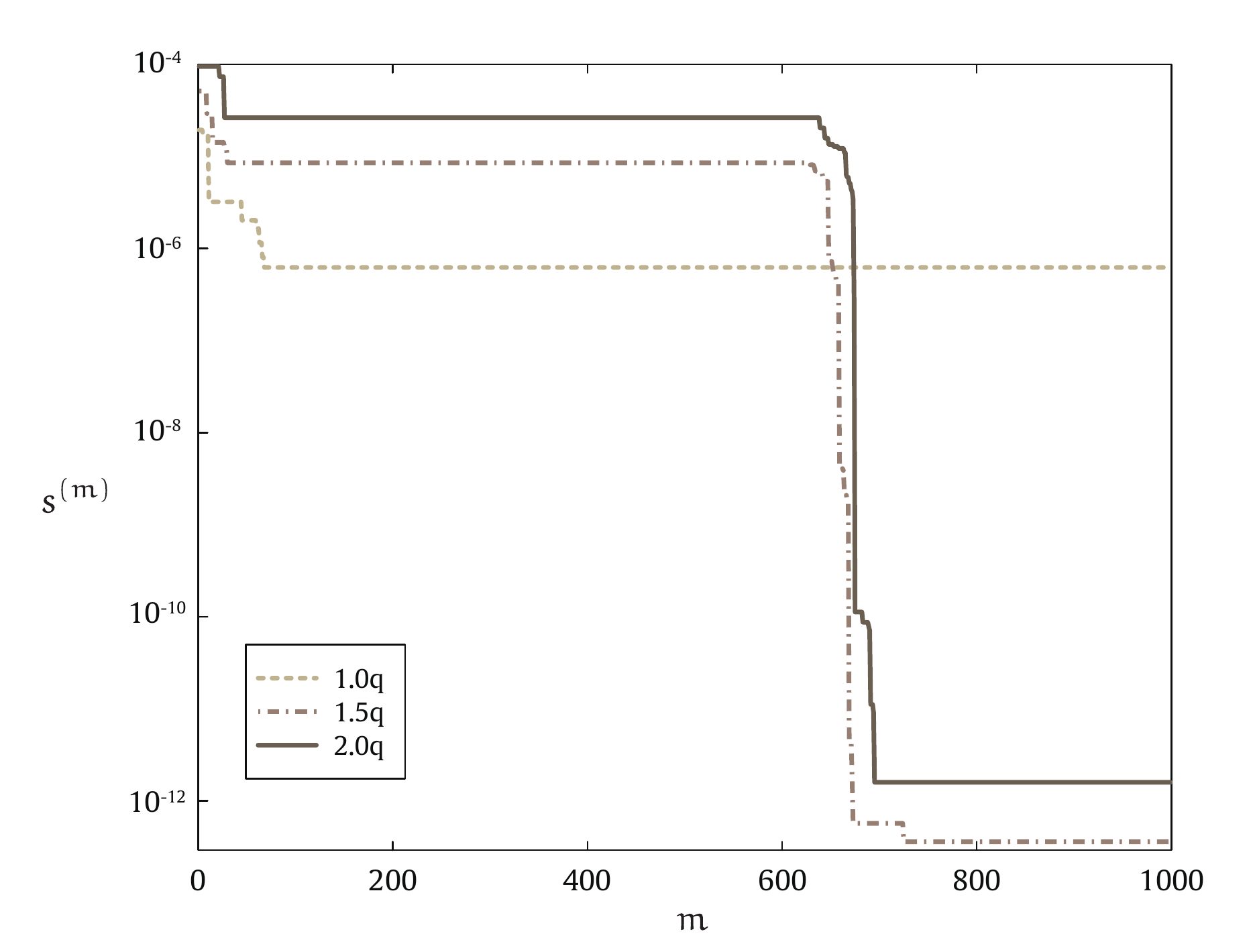}
			\caption{Step size updates in backtracking.} 
			\label{fig:steps} 
		\end{center}
	\end{figure}

	Furthermore, Fig. \ref{fig:steps} shows how the step sizes are updated by backtracking during the iterations, with three different demand levels of the Sioux Falls network as tested in the numerical experiments in Section \ref{sec:experiments}. This result demonstrates the benefit of the backtracking procedure: it allows for a fast update of the solution with large step sizes at the beginning, and later enables a stable search for the optimal solution with small step sizes.
	
	It is worth noting that, though the update rate $\xi$ is set to $0.95$ for this detailed analysis, it can be defined as smaller values in practice for efficient computation. Even if $\xi$ is small, the backtracking procedure still effectively work as it provides different values of the step size $s$ during iterations, i.e., both the fast solution search and smooth convergence. In the experiments in Sections \ref{sec:comparison} and \ref{sec:realsize}, we indeed set $\xi = 0.25$, with which AGP showed satisfying smooth convergence to the optimal solution.
	
	\section{Note on the probit assignment} \label{app:probit}
\begin{figure}[htbp]
		\begin{center}
			\includegraphics[width=10cm]{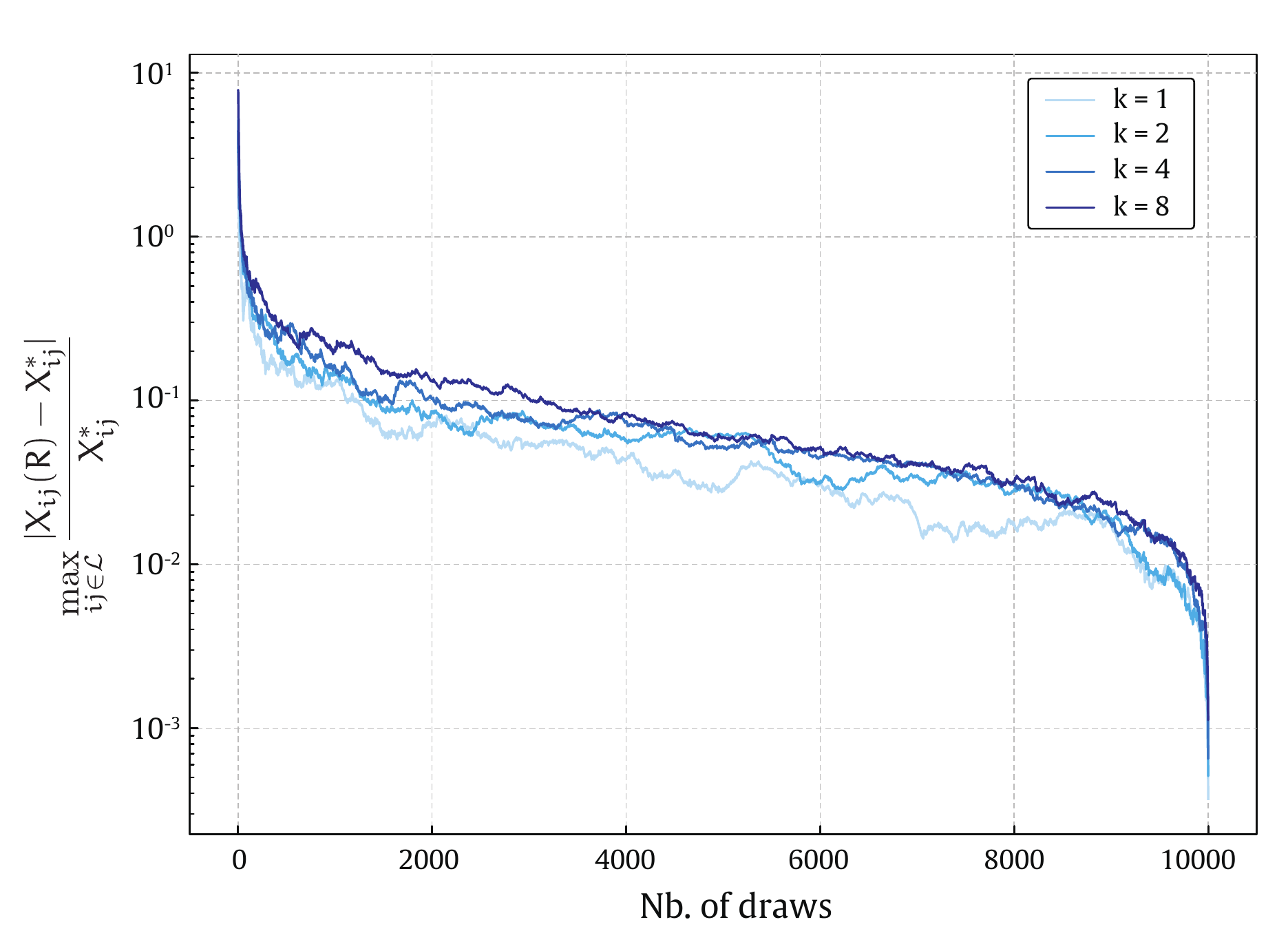}
			\caption{Approximation error of the probit assignment in different sizes of grid networks.}
			\label{fig:probit} 
		\end{center}
	\end{figure}

\noindent
	We show here an experiment to see the approximation error of the probit assignment using the grid networks of Fig. \ref{fig:largenet}. The probit assignment relies on the Monte-Carlo simulation that draws link travel costs from a multivariate normal distribution and performs the SP assignment for each draw. The expected link flow is then approximated by taking the mean $X_{ij}(R) \approx \sum_{r=1}^R X_{ij}^{(r)}/R$ where $R$ is the number of draws. In this experiment, we set the solution to $X^*_{ij} \equiv X_{ij}(10000)$ and see how the approximation error $\max_{ij \in \mathcal{L}} [ |X_{ij}(R) - X^*_{ij}|/X^*_{ij} ] $ decreases as $R$ grows. Fig. \ref{fig:probit} reports the results of four different sizes of grid networks. These results show that more than $1000$ draws are required to achieve an error within $10$\% for any network size. As expected, larger-scale networks incur more significant approximation errors. That is because there are many feasible paths between an OD pair in such networks. This experiment shows that, even if the probit model provides a more flexible description of the path correlation, the advantage may vanish when the number of draws is limited for efficiency.
	Moreover, even the probit assignment with a limited number of draws is not comparable to assignment models with a closed-form expression in terms of required computational effort. Dial's algorithm performs the SP assignment only once, and the MTA computation requires a similar effort to Dial's algorithm. An additional draw of the probit assignment requires the same level of computational effort as an assignment model with a closed-form expression. \ref{app:loadingtime} provides a comparison of computational time of loading algorithms.
	
\begin{figure}[htbp]
		\begin{center}
			\includegraphics[width=10cm]{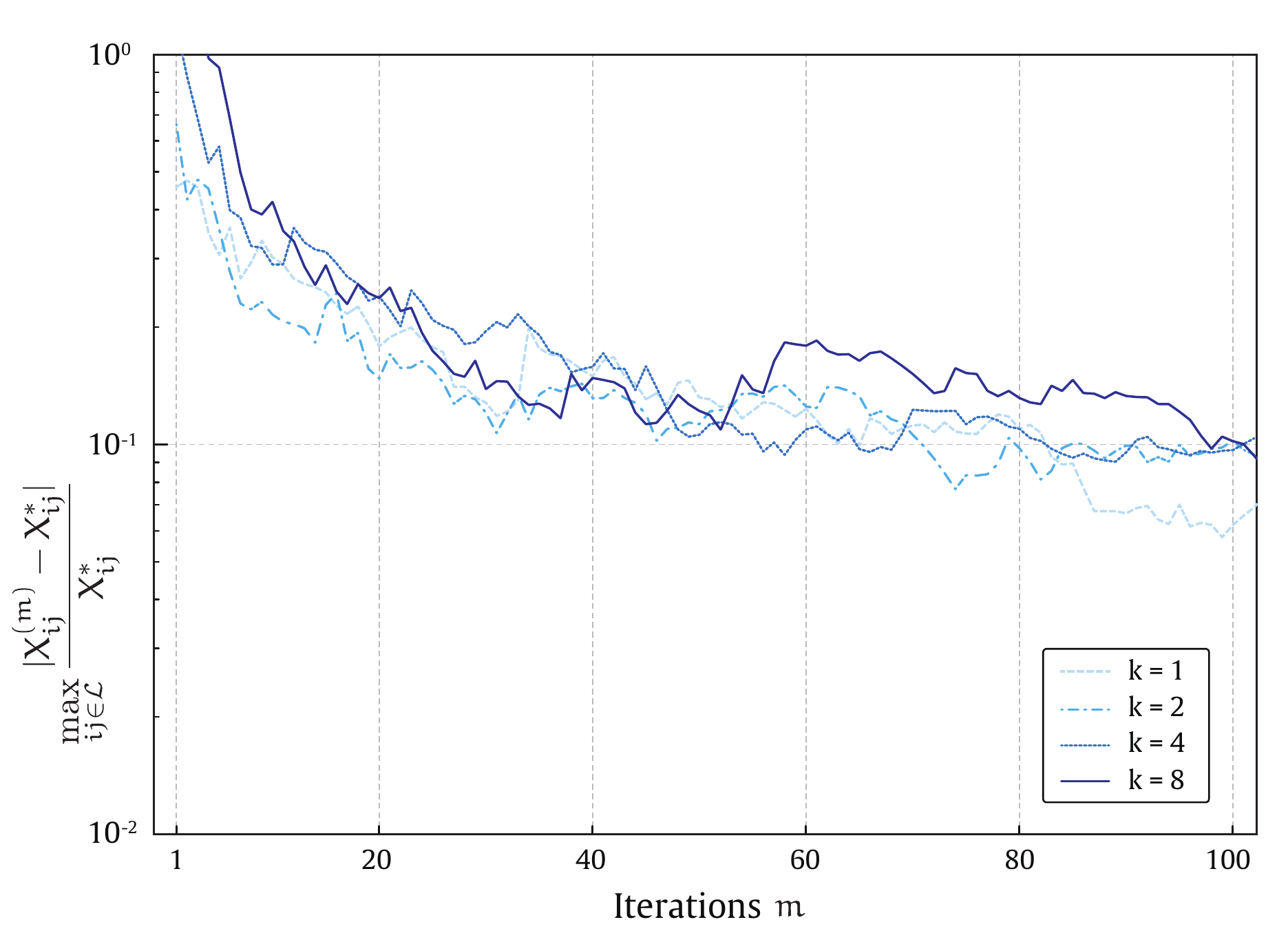}
			\caption{Convergence processes of the probit equilibrium assignment in different sizes of grid networks.}
			\label{fig:probitSUE} 
		\end{center}
	\end{figure}

	Next, we show the case of equilibrium assignment with the probit model implemented. As the reference point, we first performed the probit equilibrium assignment with $R = 10000$ draws for each iteration, solved it by MSA ($100$ iterations), and obtained $\mathbold{X}^*$. We then analyze the relative error $\max_{ij \in \mathcal{L}} [ |X^{(m)}_{ij} - X^*_{ij}|/X^*_{ij} ]$ at each iteration $m$ when implementing the probit assignment with $R = 10$, and Fig. \ref{fig:probitSUE} shows how its trajectory behaves over iterations. Even in a small network of $k=1$, the probit equilibrium assignment with $R = 10$ never achieves an error within $5$\%. In the cases of $k=2$, $4$ and $8$, the error fluctuates around $10$\%. The approximation error may become much more significant in large-scale networks where there exist almost uncountable feasible paths.

	\newpage

	\bibliographystyle{elsarticle-harv}
	\bibliography{NGEV}
	
\end{document}